\documentstyle[12pt,amsmath,amssymb,twoside]{article}

\begin{document}
\newtheorem{lem}{Lemma}[section]
\newtheorem{th}{Theorem}[section]
\newtheorem{prop}{Proposition}[section]
\newtheorem{rem}{Remark}[section]
\newtheorem{define}{Definition}[section]
\newtheorem{cor}{Corollary}[section]

\makeatletter\@addtoreset{equation}{section}\makeatother
\def\theequation{\arabic{section}.\arabic{equation}}

\newcommand{\ext}{{mk}}
\newcommand{\N}{{\Bbb N}}
\newcommand{\C}{{\Bbb C}}
\newcommand{\Z}{{\Bbb Z}}
\newcommand{\R}{{\Bbb R}}
\newcommand{\Rp}{{\R_+}}
\newcommand{\eps}{\varepsilon}
\newcommand{\Om}{\Omega}
\newcommand{\OM}{\Om_X^{\Rp}}
\newcommand{\om}{\omega}
\newcommand{\PR}{\operatorname{Pr}}
\newcommand{\XR}{{X\times\Rp}}
\newcommand{\supp}{\operatorname{supp}}
\newcommand{\la}{\langle}
\newcommand{\ra}{\rangle}
\newcommand{\Oc}{{\cal O}_{\mathrm{c}}}
\newcommand{\Bc}{{\cal B}_{\mathrm{c}}}
\newcommand{\Lext}{\Lambda_{{\mathrm \ext}}}
\newcommand{\tsigma}{{\tilde \sigma}}
\renewcommand{\emptyset}{\varnothing}
\renewcommand{\tilde}{\widetilde}
\newcommand{\Got}{{\frak G}}
\newcommand{\got}{{\frak g}}
\newcommand{\Diff}{\operatorname{Diff}_0(X)}
\newcommand{\RtX}{\R_+^X}
\newcommand{\rom}[1]{{\rm #1}}
\newcommand{\FC}{{\cal F}C_{\mathrm b}^\infty\big({\frak D},\OM\big)}
\newcommand{\vph}{\varphi}
\newcommand{\fii}{\vph}
\newcommand{\di}{\partial}
\renewcommand{\div}{\operatorname{div}}
\newcommand{\Bhat}{\widehat{{\cal B}}}
\newcommand{\un}{{\underline n}}
\newcommand{\pit}{{\pi_\tsigma}}
\newcommand{\FP}{{\cal F}{\cal P}({\frak D},\OM)}
\newcommand{\FCp}{{\cal F}C^\infty_{\mathrm p}({\frak D},\OM)}
\newcommand{\EOM}{{\cal E}^\Omega_{\pi_\tsigma}}
\newcommand{\HOM}{H^\Om_{\pi_\tsigma}}
\newcommand{\HXR}{H_\tsigma^\XR}

\thispagestyle{empty}

$\mbox{}$\vspace{10mm}
\begin{center}
\LARGE\bf Analysis and geometry\\[2mm]
\LARGE\bf on $\R_+$-marked configuration space\\[20mm]
\large\bf Yuri G. Kondratiev\\[3mm]
\large Inst.\ f.\ Angew.\ Math., Univ.\ Bonn, 53115 Bonn, Germany;\\
\large BiBoS, Univ.\ Bielefeld, 33615 Bielefeld, Germany; and\\
Inst.\ Math., NASU, 252601 Kiev, Ukraine\\[15mm]
\large\bf Eugene W. Lytvynov\\[3mm]
\large BiBoS, Univ.\ Bielefeld, 33615 Bielefeld, Germany\\[15mm]
\large\bf Georgi F. Us\\[3mm]
\large Dep.\ of Mech.\ and Math., Kiev Univ., 252033 Kiev, Ukraine
\vspace{30mm}
\end{center}
\setcounter{section}{-1}
\setcounter{page}{0}

\newpage

\begin{center}\bf Abstract\end{center}
\noindent\begin{small}

We carry out analysis and geometry on a marked configuration space $\OM$ over a
Riemannian manifold $X$ with marks from the space $\Rp$ as a natural
generalization of the work {\bf [}{\it J. Func.\ Anal}.~{\bf 154} (1998),
  444--500{\bf ]}. As a transformation group $\Got$ on $\OM$  we take the ``lifting''
to $\OM$ of the action on $\XR$ of the semidirect product of the group $\Diff$
of diffeomorphisms on $X$ with compact support and the group $\R_+^X$ of smooth
currents, i.e., all $C^\infty$ mappings of $X$ into $\Rp$ which are equal to one
outside a compact set. The marked Poisson measure $\pi_\tsigma$ on $\OM$ with
L\'evy measure $\tsigma$ on $\XR$ is proven to be quasiinvariant under the
action of $\Got$. Then, we derive a geometry on $\OM$ by a natural ``lifting''
of the corresponding geometry on $\XR$. In particular, we construct a
gradient $\nabla^\Omega$ and divergence $\div^\Om$. The associated
volume elements, i.e., all probability measures $\mu$ on $\OM$ with respect to
which $\nabla^\Om$ and $\div^\Om$ become dual operators on $L^2(\OM;\mu)$ are
identified as the mixed Poisson measures with mean measure equal to a multiple
of $\tsigma$. As a direct consequence of our results, we obtain  marked Poisson
space representations of the group $\Got$ and its Lie algebra $\got$. We
investigate also Dirichlet forms and Dirichlet operators connected with (mixed)
marked Poisson measures. In particular, we obtain conditions of ergodicity of
the semigroups generated by the Dirichlet operators. A possible generalization
of the results of the paper to the case where the marks belong to a homogeneous
space of a Lie group is noted.

\end{small}
\newpage

\tableofcontents
\newpage

\section{Introduction}
In recent few years, stochastic analysis and differential geometry on
configuration spaces have been considerably developed in a series of papers by
S.~Albeverio et~al.\ [2--5], see also \cite{Roe}. It has been shown that the
geometry of the configuration space $\Gamma_X$ over a Riemannian manifold $X$
can be constructed via a simple ``lifting procedure'' and is completely
determined by the Riemanian structure of $X$. The mixed Poisson measures
appeared to be exactly ``volume elements'' corresponding to the differential
geometry introduced on $\Gamma_X$. Intrinsic Dirichlet forms an operators, their
canonical processes, Gibbs measures on configuration spaces, integration by
parts characterization of canonical Gibbs measures, stochastic dynamics
corresponding to Gibbs measures as well as many other problems were treated in
the above framework.

A starting point for this analysis, more exactly, for the defintion of
differentiation on the configuration space, were the representation of the group
of diffeomorphisms $\Diff$ on $X$ with compact support that was constructed by
A.~M.~Vershik et~al.\ \cite{VGG2} (see also \cite{Ob,Sh,Ismagilov}) and the fact,
following from the Skorokhod theorem, that the Poisson measure is quasiinvariant
with respect to the group $\Diff$.

On the other hand, starting with the same work \cite{VGG2}, many researchers
consider representations also on marked (or compound) Poisson spaces. At the
same time, in statistical mechanics of continuous systems, marked Poisson
measures and their Gibbsian perturbations are used for the description of many
concrete models, see e.g.\ \cite{AGL78}. Hence, it is natural to ask about
geometry and analysis on these spaces. The first work in this direction was the
paper \cite{KSS}, in which, just as in the case of the usual Poisson measure,
the action of the group $\Diff$ was used for the definition of the
differentiation. However, this group proved to be too small for reconstructing
mixed marked Poisson measures as ``volume elements,'' which means that $\Diff$
is to be extended in a proper way.

Let us recall that the configuration space $\Gamma_X$ is defined as
\[ \Gamma_X:=\big\{\,\gamma\subset X\mid \#(\gamma\cap K)<\infty\text{ for each
  compact $K\subset X$}\,\big\},\]
where $\#(\cdot)$ denotes the cardinality of a set. Then, the marked
configuration space $\Omega_X^M$ over $X$ with marks from, generally speaking, a
manifold $M$ is defined as
\[ \Omega_X^M:=\big\{\,(\gamma,m)\mid \gamma\in \Gamma_X,\, m\in M^\gamma\,\big\},\]
where     $M^\gamma$ stands for the set of all maps $\gamma\ni x\mapsto m_x\in
M$. Let $\tsigma$ be a Radon measure on $ X\times M$ such that $\tsigma(K\times
M)<\infty$ for each compact $K\subset X$ and $\tsigma$ is nonatomic in $X$,
i.e., $\tsigma(\{x\}\times M)=0$ for each $x\in X$. Then, one can define on
$\Omega_X^M$ a marked Poisson measure $\pi_\tsigma$ with L\'evy measure
$\tsigma$.

Of course, one could consider $\pi_\tsigma$ as a usual Poisson measure on  the
configuration space $\Gamma_{X\times M}$ over the Cartesian product of the
underlying manifold $X$ and the space of marks $M$, and study the properties of
this measure using the results of [2--5]. However, such an approach does not
distinguish between the two different natures of $X$ and $M$ and the different
roles that these play in physics. Thus, our aim is to introduce and study such
transformations of the marked configuration space which do ``feel'' this
difference and lead to an appropriate stochastic analysis and differential
geometry.

In this paper, we will be concerned with  the model case where the space of marks $M$ is
just $\R_+$. In our forthcoming papers \cite{KLU1,KLU2}, we will generalize our results
to the case where $M$ is a homogeneous space of a Lie group and the marked
Poisson measures are replaced with Gibbs measures of Ruelle-type on $\OM$
(compare with \cite{AKR2,Kuna}).

Let $\R_+^X$ denote the group of smooth currents, i.e., all $C^\infty$ mappings
$X\ni x\mapsto \theta(x)\in\Rp$ which are equal to one outside a compact set
(depending on $\theta$). We define the group $\Got$  as the semidirect  product
of the groups $\Diff$ and $\R_+^X$: for $g_1=(\psi_1,\theta_1)$ and $g_2=(\psi_2,
\theta_2)$, where $\psi_1,\psi_2\in \Diff$ and $\theta_1,\theta_2\in\R_+^X$, the
multiplication of $g_1$ and $g_2$ is given by
\[ g_1g_2=(\psi_1\circ\psi_2,\theta_1(\theta_2\circ\psi_1^{-1})).\]
The group $\Got$ acts in $\XR$ as follows: for any $g=(\psi, \theta)\in\Got$
\[ \XR\ni(x,s)\mapsto g(x,s)=(\psi(x),\theta(\psi(x))s)\in\XR.\]
Since each $\om\in\OM$ can be interpreted as a subset of $\XR$, the action of
$\Got$ can be lifted to an action on $\OM$. The marked Poisson measure
$\pi_\tsigma$ is proven to be quasiinvariant under it. Thus, we can easily construct, in
particular, a representation of $\Got$ in $L^2(\pi_\tsigma)$.

We note that the groups of smooth (as well as measurable and continuous)
currents are classical objects in representation theory, see e.g.\
\cite{VGG1,GGV1,GGV2,VGG3,Ismagilov} and references therein for different
representations of these groups. It should be stressed, however, that our
representation of $\Got$ is reducible, because so is the regular representation
of $\Got$ in $L^2(\tsigma)$ (see subsec.~3.5 for details).

Thus, having introduced the action of the group $\Got$ on $\OM$, we proceed to
derive analysis and geometry on $\OM$ in a way parallel to the work \cite{AKR},
dealing with the usual configuration space $\Gamma_X$. In particular, we note
that the Lie algebra $\got$ of the group $\Got$ is given by $\got=V_0(X)\times
C_0^\infty(X)$, where $V_0(X)$ is the algebra of $C^\infty$ vector fields on $X$
having compact support and $C_0^\infty(X)$ is the algebra of $C^\infty$ functions
from $X$ into $\R$  with compact support. For each $(v,a)\in\got$, we define the
notion of a directional derivative of a function $F\colon\OM\to\R$ along
$(v,a)$, which is denoted by $\nabla^\Omega_{(v,a)}F$. We obtain an explicit
form of this derivative on a special set $\FC$ of smooth cylinder functions on
$\OM$, which, in turn, motivates our definition of a tangent bundle $T(\OM)$ of
$\OM$, and of a gradient $\nabla^\Omega F$. We
note only that the tangent space $T_\om(\OM)$ to the marked configuration space
$\OM$ at a point $\om=(\gamma,s)\in\OM$ is given by
\[ T_\om(\OM):= L^2(X\to T(X)\dotplus \R;\gamma).\]

Next, we derive an integration by parts formula on $\OM$, that is, we get an
explicit formula for the dual operator $\div^\Om$ of the gradient $\nabla^\Om$
on $\OM$. We prove that the probability measures on $\OM$ for which $\nabla^\Om$
and $\div^\Om$ become dual operators (with respect to
$\la\cdot,\cdot\rangle_{T(\OM)}$) are exactly the mixed marked Poisson measures
\[\mu_{\nu,\tsigma}=\int_\Rp
\pi_{z\tsigma}\,\nu(dz),\]
where $\nu$ is a probability measure on $\Rp$ (with finite first moment) and
$\pi_{z\tsigma}$ is the marked Poisson measure on $\OM$ with L\'evy measure
$z\tsigma$, $z\ge0$. This means that the mixed marked Poisson measures are
exactly the ``volume elements'' corresponding to our differential geometry on
$\OM$.

Thus, having identified the right volume elements on $\OM$, we introduce for
each measure $\mu_{\nu,\tsigma}$ the first order Sobolev space
$H_0^{1,2}(\OM,\mu_{\nu,\tsigma})$ by closing the corresponding Dirichlet form
\[ {\cal E}^\Om_{\mu_{\nu,\tsigma}}(F,G)=\int_{\OM}\la \nabla^\Om F,\nabla^\Om
G\ra_{T(\OM)}\, d\pi_{\nu,\tsigma},\qquad F,G\in\FC,\]
on $L^2(\OM,\mu_{\nu,\tsigma})$. Just as in the analysis on the usual
configuration space, this is the step where we really start doing  generic
infinite dimensional analysis. The corresponding
Dirichlet operator is denoted by $H_{\mu_{\nu,\tsigma}}^\Om$; it is a positive
definite selfadjoint operator on $L^2(\OM,\mu_{\nu,\tsigma})$. The heat semigroup $\big(
\exp(-t H^\Om_{\mu_{\nu,\tsigma}})
\big)_{t\ge0}$  generated by it is calculated explicitly. The results on the
ergodicity of this semigroup are absolutely analogous to the corresponding
results of \cite{AKR}. Particularly, we have ergodicity if and only if
$\mu_{\nu,\tsigma}=\pi_{z\tsigma}$ for some $z>0$, i.e., $\mu_{\nu,\tsigma}$ is
a (pure) marked Poisson measure.

We also clarify the relation between the intrinsic geometry on $\OM$ we
have constructed with another kind of extrinsic geometry on $\OM$ which is based
on fixing the marked Poisson measure $\pi_\tsigma$ and considering the unitary
isomorphism between $L^2(\OM,\pi_\tsigma)$ and the corresponding Fock space
\[ {\cal F}(L^2(\XR;\tsigma))=\bigoplus_{n=0}^\infty\hat
L^2((\XR)^n,n!\,\tsigma^{\otimes n}),\]
where $\hat L^2((\XR)^n,n!\,\tsigma^{\otimes n})$ is the subspace of symmetric
functions from \linebreak $L^2((\XR)^n,n!\,\tsigma^{\otimes n})$. Our main result here is to
prove that $H^\Om_{\pi_\tsigma}$ is unitarily equivalent (under the above
isomorphism) to the second quantization operator of the Dirichlet operator
$H_{\tsigma}^{\XR}$ on the
$L^2(\XR;\tsigma)$
 space.

As a consequence of the results of this paper, we obtain a representation on the
marked Poisson space $L^2(\pi_\tsigma)$ not only of the group $\Got$, but also
of its Lie algebra $\got$.

Finally, we note that one can construct in a natural way a bijection between
the marked configuration space $\OM$ and the ``compound configuration'' space
$\Om_X$ over $X$ defined as
\[ \Om_X=\{\upsilon(\cdot)=\sum_{x\in\gamma}s_x\eps_x(\cdot)\mid
\gamma\in\Gamma_X,\, s_x\in\Rp\}.\]
The image of the marked Poisson measure under this bijection is the compound
Poisson measure on $\Om_X$. So, all results of this paper can be easily
reformulated in terms of compound Poisson measures.

\section{Marked Poisson measures}

\subsection{Marked configuration space}

Let $X$ be a connected, oriented $C^\infty$ (non-compact) Riemannian manifold. The
configuration space $\Gamma_X$ over $X$ is defined as the set of all locally
finite subsets in $X$:
\[\Gamma_X:=\big\{\,\gamma\subset X\mid\#(\gamma\cap K)<\infty\ \text{for each
  compact } K\subset X\,\big\},\]
where $\#(\cdot)$ denotes the cardinality of a set. One can identify any
$\gamma\in\Gamma_X$ with the positive integer-valued Radon measure
\[\sum_{x\in\gamma}\eps_x\in{\cal M}_{\mathrm p}(X)\subset{\cal M}(X),\]
where $\sum_{x\in\emptyset}\eps_x:=\text{zero measure}$ and ${\cal M}(X)$
(resp.\ ${\cal M}_{\mathrm p}(X)$) denotes the set of all positive (resp.\ positive
integer-valued) Radon measures on ${\cal B}(X)$.

Let $\Rp:=(0,+\infty)$. The marked configuration space $\OM$ over $X$ with marks
from $\Rp$ is defined as
\[ \OM:=\big\{\,\om=(\gamma,s)\mid \gamma\in\Gamma_X,\ s\in \R_+^\gamma\,\big\},\]
where $\R_+^\gamma$ stands for the set of all maps $\gamma\ni x\mapsto
s_x\in\Rp$. Equivalently, we can define $\OM$ as
the collection of locally finite subsets in $\XR$ having the following properties:
\[\OM=\left\{\om\subset X\times\Rp\,\left|\, \begin{matrix}
 &\text{a)}\,\forall (x,s),(x',s')\in\om:(x,s)\ne(x',s')\Rightarrow x\ne x'\\
&\text{b)}\,\PR_X\om\in\Gamma_X\end{matrix}\right.\right\},\]
where $\PR_X$ denotes the projection of the Cartesian product of $X$
and $\Rp$ onto $X$. Again, each $\om \in\OM$ can be identified with the measure
\[\sum_{(x,s)\in\om}\eps_{(x,s)}\in{\cal M}_{\mathrm p}(X\times\Rp)\subset{\cal M}(X\times\Rp).\]

It is worth noting that, for any bijection $\phi\colon \XR\to\XR$, the image of
the measure $\om(\cdot)$ under the mapping $\phi$, $(\phi^*\om)(\cdot)$,
coincides with $(\phi(\om))(\cdot)$, i.e.,
\[ (\phi^*\om)(\cdot)=(\phi(\om))(\cdot),\qquad \om\in\OM,\]
where $\phi(\om)=\{\phi(x,s)\mid(x,s)\in\om\}$ is the image of $\om$ as a
subset of $\XR$.

Let $\Bc(X)$ and $\Oc(X)$ denote the families of all Borel, resp.\ open subsets
of $X$ that have compact closure. Let also $\Bc(\XR)$ denote the family of all
Borel subsets of $\XR$ whose projection on $X$ belongs to $\Bc(X)$.

Denote by $C_{\text{0,b}}(\XR)$ the set of real-valued bounded continuous
functions $f$ on $\XR$ such that
$\supp f\in\Bc(\XR)$.
 As usually, we set for any $f\in C_{\text{0,b}}(\XR)$ and $\om\in\OM$
\[\la f,\om\ra=\int_\XR f(x,s)\,\om(dx,ds)=\sum_{(x,s)\in\om}f(x,s).\]
Notice that, because of the definition of $\OM$, there are  only a finite number
of addends in the latter series.

Now, we are going to discuss the measurable structure of the space $\OM$. We
will use a ``localized'' description of the Borel $\sigma$-algebra ${\cal
  B}(\OM)$ over $\OM$.

For $\Lambda\in\Oc(X)$, define
\[\Om_\Lambda^{\Rp}:=\big\{\,\om\in\OM\mid\PR_X\om\subset\Lambda\,\big\}\]
and for $n\in\Z_+=\{0,1,2,\dots\}$
\[ \Om_\Lambda^{\Rp}(n):=\big\{\, \om\in\Om_\Lambda^{\Rp}\mid\#(\om)=n\,\big\}.\]
It is obvious that
\[\Omega_\Lambda^{\Rp}=\bigsqcup_{n=0}^\infty \Omega_\Lambda^\Rp(n).\]

Let
\[\Lext:=\Lambda\times\Rp\]
(i.e., $\Lext$ is the set of all ``marked'' elements of  $\Lambda$)
and let
\[{\widetilde \Lambda}^n_{\text{\ext}}:=\big\{
((x_1,s_1),\dots,(x_n,s_n))\in\Lambda_{\text{\ext}}^n\mid x_j\ne x_k\ \text{if }j\ne k\,\big\}.
\]
There is a bijection
\begin{equation}\label{one}
{\cal L}_\Lambda^{(n)}\colon {\widetilde\Lambda}_{\text{\ext}}^n/{\frak
  S}_n\mapsto
\Om_{\Lambda}^{\Rp}(n)\end{equation}
given by
\[ {\cal L}_\Lambda^{(n)}\colon((x_1,s_1),\dots,(x_n,s_n))\mapsto
\{(x_1,s_1),\dots,(x_n,s_n)\}\in \Om_{\Lambda}^\Rp(n),\]
where ${\frak S}_n$ is the permutation group over $\{1,\dots,n\}$. On
$\Lext^n/{\frak S}_n$ one introduces the related metric
\begin{align*}
&
\delta\big[((x_1,s_1),\dots,(x_n,s_n)),((x_1',s_1'),\dots,(x_n',s_n'))\big]\\
&\qquad =\inf_{\pi\in
{\frak S}_n}
d^n\big[((x_1,s_1),\dots,(x_n,s_n)),((x'_{\pi(1)},s'_{\pi(1)}),\dots,(x'_{\pi(n)},
s'_{\pi(n)}))\big],
\end{align*}
where $d^n$ is the metric on $\Lext^n$ driven from the original metrics on $X$
and $\Rp$. Then, ${\widetilde \Lambda}_{\text{\ext}}^n/{\frak S}_n$ becomes an open set in
$\Lext^n/{\frak S}_n$ and let ${\cal B}({\widetilde\Lambda}_{\text{\ext}}^n/{\frak S}_n)$ be the
trace $\sigma$-algebra on ${\widetilde \Lambda}_{\text{\ext}}^n/{\frak S}_n$ generated by ${\cal
  B}(\Lext^n/{\frak S}_n)$. Let then ${\cal B}(\Om_\Lambda^\Rp(n))$ be the image
$\sigma$-algebra of ${\cal B}({\widetilde\Lambda}_{\text{\ext}}^n/{\frak S}_n)$ under the bijection
${\cal L}_\Lambda^{(n)}$ and let ${\cal B}(\Om_\Lambda^{\Rp})$ be the
$\sigma$-algebra on $\Om_\Lambda^{\Rp}$ generated by the usual topology of
(disjoint) union of topological spaces.

For any $\Lambda\in\Oc(X)$, there is a natural restriction map
\[ p_\Lambda\colon\OM\mapsto \Om_\Lambda^{\Rp}\]
defined by
\[\OM\ni\om\mapsto p_\Lambda (\om):=\om\cap\Lext\in\Om_\Lambda^\Rp.\]
The topology on $\OM$ is defined as the weakest topology making all the mappings
$p_\Lambda$ continuous. The associated $\sigma$-algebra is denoted by ${\cal
  B}(\OM)$.

For each $B\in\Bc(\XR)$, we introduce a function
$N_B\colon\OM\to\Z_+=\{0,1,2,\dots\}$ such that
\begin{equation} N_B(\om):=\#(\om\cap B),\qquad \om\in\OM.\label{NB}\end{equation}
Then, it is not hard to see that ${\cal B}(\OM)$ is the smallest
$\sigma$-algebra on $\OM$ such that all the functions $N_B$ are
measurable.

\subsection{Marked Poisson measure}

In order to construct a marked Poisson measure, we fix:

(i) an intensity measure $\sigma$ on the underlying manifold $X$, which is
supposed to be  a   nonatomic Radon one,

(ii) a non-negative function
\[X\times{\cal B}(\Rp)\ni (x,\Delta)\mapsto p(x,\Delta)\in\R_+
\]
such that, for $\sigma$-a.a.\ $x\in X$, $p(x,\cdot)$ is a finite Radon measure
on $\Rp$.

Now, we define a measure $\tsigma$ on $(\XR,{\cal B}(\XR))$ as follows:
\begin{equation}\label{1}
\tsigma(A)=\int_{A} p(x,ds)\,\sigma(dx),\qquad A\in{\cal B}(\XR) .
\end{equation}
We will suppose that the measure $\tsigma$ is infinite and  for any
$\Lambda\in\Bc(X)$
\begin{equation}\label{2}
\tsigma(\Lext)=\int_X {\pmb1}_\Lambda(x)p(x,\Rp)\,\sigma(dx)<\infty,\end{equation}
i.e., $p(x,\Rp)\in L^1_{\mathrm loc}(\sigma)$.

Now, we wish to introduce a marked Poisson measure on $\OM$ (cf.\ e.g.\
\cite{kingman,Kall}).  To this end, we
take first the measure $\tsigma^{\otimes n}$ on $(X\times\Rp)^n$, and for any
$\Lambda\in\Oc(X)$, $\tsigma^{\otimes n}$ can be considered as a finite measure
on $\Lambda_{\text{\ext}}^n$. Since $\sigma$ is nonatomic, we get
$\tsigma^{\otimes 2}(D_\Lambda)=0$, where
\[ D_\Lambda=\big\{\,((x_1,s_1),(x_2,s_2))\in \Lambda^2_{\text{\ext}}\mid
x_1=x_2\,\big\}=\big\{\,(x_1,x_2)\in X^2\mid x_1=x_2\,\big\}\times\R_+^2.\]
Therefore,
\[\tsigma^{\otimes n}
(\Lambda_{\text{\ext}}^n\setminus {\widetilde\Lambda}_{\text{\ext}}^n)=0
\]
and we can consider $\tsigma^{\otimes n}$ as a measure on
$(\widetilde\Lambda^n_{\text{\ext}}/ {\frak S}_n,{\cal
  B}(\widetilde\Lambda^n_{\text{\ext}}/{\frak S}_n))$
such that
$$\tsigma^{\otimes n}(\widetilde\Lambda^n_{\text{\ext}}/{\frak
S}_n)=\tsigma(\Lext)^n.$$

Denote by $\tsigma_{\Lambda,n}:=\tsigma^{\otimes n}\circ ({\cal
  L}_\Lambda^{(n)})^{-1}$ the image measure on $\Omega_\Lambda^\Rp(n)$ under the
  bijection \eqref{one}. Then, we can define a measure
  $\lambda^\Lambda_{\tsigma}$ on $\Omega_\Lambda^{\Rp}$ by
\[ \lambda_\tsigma^\Lambda:=\sum_{n=0}^\infty\frac1{n!}\,\tsigma_{\Lambda,n},\]
where $\tsigma_{\Lambda,0}:=\eps_{\emptyset}$ on
$\Omega_\Lambda^\Rp(0)=\{\emptyset\}$. The measure $\lambda_\tsigma^\Lambda$ is
finite and
$\lambda^\Lambda_\tsigma(\Omega_\Lambda^\Rp)=e^{\tsigma(\Lext)}$. Hence, the
measure
\[ \pi_\tsigma^\Lambda:=e^{-\tsigma(\Lext)}\lambda^\Lambda_{\tsigma}\]
is a probability measure on ${\cal B}(\Omega_\Lambda^\Rp)$. It is not hard to
check the consistency property of the family
$\{\pi_\tsigma^\Lambda\mid\Lambda\in\Oc(X)\}$ and thus to obtain a unique
probability measure $\pi_\tsigma$ on ${\cal B}(\OM)$ such that
\[ \pi_\tsigma^\Lambda=p_\Lambda^*\pi_\tsigma,\qquad \Lambda\in\Oc(X).\]
This measure $\pi_\tsigma$ will be called a marked Poisson measure with L\'evy
measure $\tsigma$.

For any function $\fii\in C_{\text{0,b}}(\XR)$, it is easy to calculate the Laplace
transform of the measure $\pi_\tsigma$
\begin{equation}\ell_{\pi_\tsigma}(\fii):=\int_{\OM}e^{\la
    \fii,\om\ra}\,\pi_\tsigma(d\om)
=\exp\bigg(\int_{\XR}(e^{\fii(x,s)}-1)\,\tsigma(dx,ds)\bigg).\label{4}\end{equation}
Particularly, if $C_0(X)$ denotes the space of continuous functions on $X$ with
compact support, we obtain, for each $u\in C_0(X)$ such that $u(x)\le 0$ for
$\sigma$-a.a.\ $x\in X$,
\begin{align*}
L_{\pi_\tsigma}(u):&=\int_{\OM}\exp\big [ \la su,\om\ra\big]\,\pi_{\tsigma}(d\om)\\
&=\int_{\OM}\exp\bigg[\sum_{x\in\gamma}s_xu(x)\bigg]\pi_\tsigma(d(\gamma,s))\\
&=\exp\bigg(\int_\XR(e^{su(x)}-1)\,\tsigma(dx,ds)\bigg).
\end{align*}
This formula is also sufficient to define $\pi_\tsigma$. \vspace{2mm}

{\bf Example 1.} Let $p(x,\cdot)\equiv \eps_1(\cdot)$, $x\in X$. Then,
$\tsigma=\sigma\otimes\eps_1$ and $\pi_\tsigma$ is just the Poisson measure on
$(\Gamma_X,{\cal B}(\Gamma_X))$ with intensity $\sigma$.

{\bf Example 2.} Let $p(x,\cdot)\equiv \tau(\cdot)$, $x\in X$, where $\tau$ is a
finite measure on $(\Rp,{\cal B}(\Rp))$. Now,
$\tsigma=\hat\sigma=\sigma\otimes\tau$ and $\pi_\tsigma$ coincides with the
marked Poisson measure under consideration in \cite{KSS}.
Notice that the choice of $\tsigma=\hat\sigma$ as a product measure means a
position-independent marking, while the choice of a general $\tsigma$ of the
form \eqref{1} leads to a positive-depending marking.
\vspace{2mm}

In what follows, we will suppose that the measure $\sigma$ is equivalent to the
Riemannian volume $m$ on $X$: $\sigma(dx)=\rho(x)\,m(dx)$ with $\rho>0$ $m$-a.s.,
and that for $m$-a.a.\ $x\in X$ $p(x,\cdot)$ is equivalent to the restriction of
the Lebesgue measure to $\Rp$, denoted by $\lambda$:
\[ p(x,ds)=p(x,s)\,\lambda(ds)\quad\text{with }p(x,s)>0\ \lambda\text{-a.a.}\ s\in
\Rp.\]
Thus, the measure $\tsigma$ can be written in the form
\[\tsigma(dx,ds)=\rho(x)p(x,s)\,\lambda(ds)m(dx).\]
The condition $\tsigma(\Lext)<\infty$, $\Lambda\in\Bc(X)$, implies that the
function
\[ q(x,s):=\rho(x)p(x,s)\]
satisfies
\begin{equation}\label{c1}
q^{1/2}\in L^2_{\text{loc}}(X;m)\otimes L^2(\R_+;\lambda).\end{equation}
We will
suppose additionally that the following stronger condition is fulfilled:
\begin{equation}
\big(\max\{1,s\}q(x,s)\big)^{1/2}\in L^2_{\text{loc}}(X;m)\otimes
L^2(\R_+;\lambda).\label{c2}
\end{equation}

\section{Transformations of the marked Poisson measure}
\subsection{Group of transformations of the marked configuration space}

We are looking for a natural group $\Got$ of transformations of $\OM$ such that

(i) $\pi_{\tsigma}$ is $\Got$-quasiinvariant;

(ii) $\Got$ is big enough to reconstruct $\pi_{\tsigma}$ by the Radon--Nikodym
density $\dfrac{dg^*\pi_\tsigma}{d\pi_\tsigma}$, where $g$ runs through $\Got$.

Let us recall that in the work \cite{KSS} the group $\Diff$ was taken as $\Got$,
just in the same way as in the case of the usual Poisson measure
\cite{AKR}. Here, $\Diff$ stands for the group of diffeomorphisms of $X$ with
compact support, i.e., each $\psi\in\Diff$ is a diffeomorphism of $X$ that is
equal to the identity outside a compact set (depending on $\psi$). The group
$\Diff$ satisfies (i). However, unlike the case of the Poisson measure, the
condition (ii) is not now satisfied, because, for example, in the case where
$\tsigma=\sigma\otimes\tau$, there is  no information about the measure $\tau$ that
is contained in $\dfrac{d\psi^*\pi_\tsigma}{d\pi_\tsigma}$, see
\cite{KSS}.

This is why we need a proper extension of the group $\Diff$. Let us consider the
group of {\it smooth  currents}, i.e., all $C^\infty$ mappings
\[X\ni x\mapsto\theta(x)\in\R_+,\]
which are equal to one outside a compact set (depending on $\theta$). A
multiplication $\theta_1\theta_2$ in this group is defined as the pointwise
multiplication of mappings $\theta_1$ and $\theta_2$. In representation theory
this group is denoted by $\RtX$, or $C^\infty_0(X;\Rp)$.

The group $\Diff$ acts in $\RtX$  by automorphisms: for each $\psi\in \Diff$,
\[\RtX\ni\theta\overset{\alpha}{\mapsto}\alpha(\psi)\theta:=\theta\circ\psi^{-1}\in\RtX.\]
 Thus, we can endow the Cartesian product of $\Diff$ and $\RtX$ with the
 following multiplication: for $g_1=(\psi_1,\theta_1)$, $g_2=(\psi_2,\theta_2)\in
 \Diff\times\RtX$
\[ g_1g_2=(\psi_1\circ\psi_2,\theta_1(\theta_2\circ\psi_1^{-1}))\]
and obtain a semidirect product
\[ \Diff\underset{\alpha}{\times}\RtX=:\Got\]
of the groups $\Diff$ and $\RtX$.

The group $\Got$ acts in $\XR$ in the following way: for any
$g=(\psi,\theta)\in\Got$
\begin{equation}\label{3}
\XR\ni(x,s)\mapsto g(x,s)=(\psi(x),\theta(\psi(x))s)\in\XR.\end{equation}
If id denotes the identity diffeomorphism of $X$ and ${\pmb{1}}$ is the
function identically equal to one on $X$, then we will just identify $\psi$ with
$(\psi,{\pmb{1}})$ and $\theta$ with $(\operatorname{id},\theta)$. The
action \eqref{3} of an arbitrary $g=(\psi,\theta)$ can be represented now as
\[(x,s)\mapsto g(x,s)=\theta\psi(x,s),\]
where
\begin{align*}
\psi(x,s)&=(\psi(x),s),\\
\theta(x,s)&=(x,\theta(x)s).\end{align*}

For any $g=(\psi,\theta)\in\Got$, denote $K_g:=K_\psi\cup K_\theta$, where
$K_\psi$ and $K_\theta$ are the minimal closed sets in $X$ outside of which
$\psi=\operatorname{id}$ and $\theta={\pmb{1}}$, respectively. Evidently,
$K_g\in{\cal B}_{\mathrm c}(X)$,
\[g(K_g)_{\mathrm \ext}=(K_g)_{\mathrm \ext},\]
and $g$ is the identity transformation outside $(K_g)_{\mathrm \ext}$.

Noting that
$$g^{-1}(x,s)=(\psi,\theta)^{-1}(x,s)=(\psi^{-1},\theta^{-1}\circ\psi)(x,s)
=(\psi^{-1}(x),\theta^{-1}(x)s),$$
we easily deduce the following

\begin{prop}\label{prop1} The measure
$\tsigma$ is $\Got$-quasiinvariant and for any $g=(\psi,\theta)\in\Got$
the Radon--Nikodym density is given by
\[
\left\{
\begin{aligned}
\text{}& p_g^{\tsigma}(x,s):=\frac{d(g^*\tsigma)}{d\tsigma}(x,s)=
\frac{q(\psi^{-1}(x),\theta^{-1}(x)
s)}{q(x,s)\theta(x)}\,J_m^\psi(x),\\
 \text{}&\text{\rom{if}
  }(x,s)\in\{0<q(x,s)<\infty\}\cap\{0<q(\psi^{-1}(x),\theta^{-1}(x)s)<\infty\},\\
\text{}&p_g^\tsigma(x,s)=1, \qquad \text{\rom{otherwise}},
\end{aligned}
\right.\]
where $J_m^\psi$ is the Jacobian determinant of $\psi$ \rom(w\rom.r\rom.t\rom.\
the Riemannian volume $m$\rom{).}
\end{prop}

\subsection{$\Got$-quasiinvariance of the marked Poisson measure}

Any $g\in\Got$ defines by \eqref{3} a transformation of $\XR$, and,
consequently, $g$ has the following ``lifting'' from $\XR$ to $\OM$:
\begin{equation}\label{kk} \OM\ni\om\mapsto g(\om)=\big\{\,
g(x,s)\mid (x,s)\in\om\,\big\}\in\OM.\end{equation}
(Note that, for a given $\om\in\OM$, $g(\om)$ indeed belongs to $\OM$ and
coincides with $\om$ for all but a finite number of points.)
The mapping \eqref{kk} is obviously measurable and we can define the image
$g^*\pi_\tsigma$ as usually.
The following proposition is an analog of a corresponding fact about Poisson measures.

\begin{prop}\label{prop2} For any $g\in\Got$\rom, we have
\[ g^*\pi_\tsigma=\pi_{g^*\tsigma}.\]
\end{prop}

\noindent{\it Proof}. The proof is the same as for the usual  Poisson
measure $\pi_\sigma$ with intensity $\sigma$ and $\phi\in\Diff$ (e.g.,
\cite{AKR}), one has just to calculate the Laplace transform of the measure
$g^*\pi_\tsigma$ for any $f\in C_{\text{0,b}}(\XR)$ and to use the formula
\eqref{4}.\quad$\blacksquare$

\begin{prop}\label{prop3} The marked Poisson measure $\pi_{\tsigma}$ is
  quasiinvariant w.r.t.\ the group $\Got$, and for any $g\in\Got$ we have
\begin{equation}\label{rr}
\frac{d(g^*\pi_\tsigma)}{d\pi_\tsigma}(\om)=\prod_{(x,s)\in\om}p_g^\tsigma(x,s)
.\end{equation}
\end{prop}

\noindent{\it Proof}. By Proposition~\ref{prop1}, the measures $g^*\tsigma$ and
$\tsigma$ are equivalent for each $g\in\Got$. To  apply the Skorokhod theorem on
absolute continuity of Poisson measures, we have to check the integrability
condition
\[\sqrt{p_g^\tsigma}-1\in L^2(X\times\Rp;\tsigma),\]
or the more restrictive one
\[ p_g^\tsigma-1\in L^1(X\times\Rp;\tsigma)\]
(see, e.g., \cite{Sk,Ta}). But the change of variables formula gives
\begin{align*}
\int_{\XR}|p_g^{\tsigma}(x,s)-1|\,\tsigma(dx,ds)
&=\int_{(K_g)_{\text{\ext}}}|p_g^\tsigma(x,s)-1|\,\tsigma(dx,ds)\\
&\le 2\tsigma((K_g)_{\mathrm \ext})<\infty.
\end{align*}
Hence, $\dfrac{d(g^*\pi_\tsigma)}{d\pi_\tsigma}$ is equal to the right hand side
of \eqref{rr} multiplied by the factor
\[
\exp\bigg[\int_{\XR}(1-p_g^\tsigma(x,s))\,\tsigma(dx,ds)\bigg].
\]
Finally, noting that
\begin{align*}
\int_{\XR}(1-p_g^{\tsigma}(x,s))\,\tsigma(dx,ds)
&=\int_{(K_g)_{\text{\ext}}}(1-p_g^\tsigma(x,s))\,\tsigma(dx,ds)\\
&= \tsigma((K_g)_{\mathrm \ext})-\tsigma((K_g)_{\mathrm \ext})=0,
\end{align*}
we conclude the proposition.\quad$\blacksquare$

\begin{rem}\rom{
Notice that only a finite (depending on $\omega$) number of factors in the product on the right hand side of \eqref{rr}
are not equal to one.
}\end{rem}

\subsection{Transformations of compound Poisson space and quasiinvariance of
  compound Poisson measures}

Let us recall that the ``compound configuration'' space $\Om_X$ over $X$ is
defined as follows (see e.g.\ \cite{KSSU,KSS}):
\begin{equation}\label{z2.4}
\Omega_X:=\{\upsilon(\cdot)=\sum_{x\in\gamma}s_x\eps_x(\cdot)\mid
\gamma\in\Gamma_X,\, s_x\in\R_+\}.\end{equation}
Notice that, for any $\upsilon\in\Omega_X$,
\begin{equation}\label{z2.5}
\{(x,s_x)\mid x\in\gamma\}=(\gamma,s)=\om\in\OM,\end{equation}
where the function $s\in\R_+^\gamma$ is defined by $s(x)=s_x$, $x\in\gamma$.

Any $\upsilon\in\Om_X$ is a Radon measure on $(X,{\cal B}(X))$. Hence, $\Om_X$
can be equipped with the Borel measurable structure as a subset of ${\cal M}(X)$
with the vague topology. Moreover, for any $u\in C_0(X)$, we have
\[\la u,\upsilon\ra=\sum_{x\in\gamma}s_x u(x)=\la su,\om\ra\]
with $\om$ of the form \eqref{z2.5}. Thus, the measurable space $(\Om_X,{\cal
  B}(\Om_X))$ is nothing but the image of the marked configuration space
  $(\OM,{\cal B}(\OM))$ under the transformation
\begin{equation}\label{z2.6}
\OM\ni\om\mapsto ({\cal
  I}\om)(\cdot)=\sum_{(x,s_x)\in\om}s_x\eps_x(\cdot)\in\Om_X,\end{equation}
which determines a one-to-one correspondence between $\OM$ and $\Om_X$. This is
  why the group $\Got$ generates the group $\widehat{\Got}=\{\hat g={\cal
  I}g{\cal I}^{-1}\mid g\in\Got\}$ of measurable transformations of the
  ``compound configuration'' space $\Om_X$.

For our fixed measure $\tsigma$ on $\XR$, we define a compound Poisson measure
$\mu_{\tsigma}^{\mathrm CP}$ as the image of the marked Poisson measure
$\pi_\tsigma$ under the transformation \eqref{z2.6}, i.e.,
\begin{equation}\label{z2.7}\mu_\tsigma^{\mathrm CP}:={\cal
    I}^*\pi_\tsigma.\end{equation}
This definition generalizes evidently the corresponding one from \cite{KSSU,KSS}
(compare with Example~2; in particular, $\mu^{\mathrm CP}_{\sigma\times \tau}$
coincides with the $\pi_\sigma^\tau$ from \cite{KSS}).

Thus, by using the bijection \eqref{z2.6} and the definition \eqref{z2.7}, all
results of this paper can be easily reformulated in terms of the compound
Poisson space $(\Om_X,{\cal B}(\Om_X),\mu_\tsigma^{\mathrm CP})$. For example,
we deduce from Proposition~\ref{prop3} the $\widehat{\Got}$-quasiinvariance of
the compound Poisson measure. Namely, the following statement holds.

\begin{prop}\label{propz2.4}
The compound Poisson measure $\mu_\tsigma^{\mathrm CP}$ is quasiinvariant with
respect to the group $\widehat {\Got}$\rom, and for any $\hat
g\in\widehat{\Got}$
\begin{equation}\label{z2.8}
\frac{d(\hat g^*\mu_\tsigma^{\mathrm CP})}{d\mu_\tsigma^{\mathrm
    CP}}(\upsilon)=\prod_{(x,s)\in
{\cal I}^{-1}\upsilon}p^\tsigma_g(x,s).\end{equation}\end{prop}

\noindent {\it Proof}. By Proposition~\ref{prop3}, for any $g\in\Got$,
$g^*\pi_\tsigma\sim\pi_\tsigma$ and the Radon--Nikodym
derivative has the
form   \eqref{rr}. Since ${\cal I}$ is a bijection, the equality \eqref{z2.7}
yields that, for any $\hat g\in\widehat {\Got}$, $\hat g^*\mu_\tsigma^{\mathrm
  CP}\sim
\mu_\tsigma^{\mathrm CP}$, and the value of $\dfrac{d\hat
  g^*\mu_\tsigma^{\mathrm CP}}{d\mu_\tsigma^{\mathrm CP}}$ at a point
$\upsilon\in\Om_X$ coincides with the value of
$\dfrac{dg^*\pi_\tsigma}{d\pi_\tsigma}$ at the point $\om={\cal
  I}^{-1}\upsilon$.\quad $\blacksquare$

\begin{rem}\rom{
Setting in \eqref{z2.8} $\tsigma=\sigma\times\tau$ and $\hat g^*={\cal
  I}(\psi,\pmb{1}){\cal I}^{-1}$, one obtains as a consequence of
  Proposition~\ref{propz2.4} the result of Proposition~2.8 from \cite{KSS}.}\end{rem}

\section{The differential geometry of marked configuration space}

\subsection{The tangent bundle of $\Omega^{\Rp}_{X}$}

Let us denote by $V_0(X)$ the set of $C^\infty$ vector fields on $X$ (i.e.,
smooth sections of $T(X)$) that have compact support. Let also
$C^\infty_0(X)$ stand for the set of all $C^\infty$ functions from $X$ into $\R$ that have
compact support. Then, $\got:=V_0(X)\times C^\infty_0(X)$ can be thought of as a Lie
algebra that corresponds to the Lie group $\Got$. More precisely, for any fixed
$v\in V_0(X)$ and for any $x\in X$, the curve
\[\R\ni t\mapsto\psi_t^v(x)\in X\]
is defined as the solution of the following Cauchy problem
\begin{equation}\label{beta}\left\{
\begin{aligned}
\frac{d}{dt}\psi_t^v(x)&=v(\psi_t^v(x)),\\
\psi_0^v(x)&=x.\end{aligned}\right.
\end{equation}
Then, the mappings $\{\psi_t^v,\, t\in\R\}$ form a one-parameter subgroup of
diffeomorphisms in $\Diff$ (see, e.g., \cite{Bo}):
\begin{align*}
1)&\forall t\in\R\quad \psi_t^v\in\Diff,\\
2)&\forall t_1,t_2\in\R\quad \psi_{t_1}^v\circ\psi_{t_2}^v=\psi^v_{t_1+t_2}.
\end{align*}
Next, for each function $a\in C^\infty_0(X)$, we define
\[ \theta_t^a=\theta_t^a(x):=e^{ta(x)}\in\R_+^X,\]
and $\{\theta_t^a,\, t\in\R\}$ form a one parameter subgroup of $\R_+^X$.

Thus, we can consider, for an arbitrary $(v,a)\in\got$, the curve
$\{(\psi_t^v,\theta_t^a),\,t\in\R\}$ in $\Got$. (Notice that this curve does not
form a subgroup of $\Got$.) Hence, to any $\om\in\OM$ there corresponds the
following curve in $\OM$:
\[\R\ni t\mapsto (\psi_t^v,\theta_t^a)\om\in\OM.\]
Define now, for a function $F\colon \OM\to\R$, the directional derivative of $F$
along $(v,a)$ as
\[ (\nabla^\Om_{(v,a)}F)(\om):=\frac
d{dt}F((\psi_t^v,\theta_t^a)\om)\big|_{t=0},\]
provided the right hand side exists. We will also denote by $\nabla^\Om_v$ and
$\nabla^\Om_a$ the directional derivatives along $(v,0)$ and $(0,a)$,
respectively.

Absolutely analogously, one defines for a function $\fii\colon \XR\to\R$ the
directional derivative of $\fii$ along $(v,a)$:
\begin{equation}\label{small}
(\nabla^\XR_{(v,a)}\fii)(x,s)=\frac{d}{dt}\fii((\psi_t^v,\theta_t^a)(x,s))\big|_{t=0}.\end{equation}
Then, for a continuously differentiable function $\fii$, we have from
\eqref{3}, \eqref{beta}, and \eqref{small}
\begin{gather}
(\nabla^\XR_{(v,a)}\fii)(x,s)=\la\nabla^X\fii(\psi_t^v(x),\theta_t^a(\psi_t^v(x))s),\frac
d{dt}\psi_t^v(x)\ra_{T_{\psi_t^v(x)}(X)}\big|_{t=0}\notag\\
\mbox{}+\frac{\di}{\di
  s}\fii(\psi_t^v(x),\theta_t^a(\psi_t^v(x))s)\,\frac d{dt}e^{ta(\psi_t^v(x))}s\big|_{t=0} \notag\\
=\la\nabla^X\fii(x,s),v(x)\ra_{T_x(X)}+s\frac\di{\di
  s}\fii(\psi_t^v(x),\theta_t^a(\psi_t^v(x))s) e^{ta(\psi_t^v(x))}\times\notag\\
\times(a(\psi_t^v(x))+t\la\nabla^X
a(\psi_t^v(x)),v(\psi_t^v(x))\ra_{T_{\psi_t^v(x)}(X)})\big|_{t=0}\notag\\
=\la\nabla^X\fii(x,s),v(x)\ra_{T_x(X)}+s\frac\di{\di s}\fii(x,s)a(x)\notag\\
=\la\nabla^\XR\fii(x,s),(v(x),a(x))\ra_{T_{(x,s)}(\XR)}.\label{9}
\end{gather}
Here,
\[ T_{(x,s)}(\XR):= T_x(X)\dotplus\R\]
and
\[\nabla^\XR:=(\nabla^X,\nabla^{\Rp}),\]
where $\nabla^X$
 denotes the gradient on $X$ and
\[\nabla^{\Rp}=s\frac \di{\di s}.\]

Let us introduce a special class of ``nice functions'' on $\OM$.
 We will say that a function $\fii\colon\XR\to\R$ is a $C^\infty$-function on
 $\XR$ if it
 is $C^\infty$ with respect to the gradient $\nabla^X$ and the usual derivative
 in $s$. Denote by $\frak D$ the set of all  $C^\infty$-functions $\fii$ on
 $\XR$ with support from $\Bc(\XR)$ such that $\fii$ and all its derivatives are
 bounded and moreover the following estimate holds for all $n\in\N$
\begin{equation}\label{trinidad}
\Big|\frac{\di^n}{\di s^n}\fii(x,s)\Big|\le C_{n,\fii}(\max\{1,s\})^{-n}.\end{equation}
Next, let
$C_{\text{b}}^\infty(\R^N)$ stand for the space of all
$C^\infty$-functions on $\R^N$ which together with all their derivatives are
bounded. Then, we can introduce $\FC$
as the set of all functions $F\colon\OM\mapsto\R$ of the form
\begin{equation}\label{5}
F(\om)=g_F(\la\vph_1,\om\ra,\dots,\la\vph_N,\om\ra),\qquad
\om\in\OM,\end{equation}
where $\vph_1,\dots,\fii_N\in{\frak D}$ and $g_F\in C_{\mathrm b}^\infty(\R^N)$
(compare with \cite{AKR}). $\FC$ will be called the set of smooth cylinder
functions on $\OM$.

For any  $F\in\FC$ of the form \eqref{5} and a  given $(v,a)\in \got$, we have,
just as in \cite{AKR},
\begin{align*}F((\psi_t^v,\theta_t^a)\om)&=g_F(\la\fii_1,(\psi_t^v,\theta_t^a)\om\ra,\dots,
\la\fii_N,(\psi_t^v,\theta_t^a)\om\ra)\\
&=g_F(\la\fii_1\circ(\psi_t^v,\theta_t^a),\om\ra,\dots,\la\fii_N\circ(\psi_t^v,\theta_t^a),\om\ra),\end{align*}
and therefore
\begin{equation}\label{8}
(\nabla_{(v,a)}^\Omega F)(\om)=\sum_{j=1}^N\frac{\di g_F}{\di r_j}(\la\fii_1,\om\ra,\dots,\la\fii_N,\om\ra)
\la\nabla_{(v,a)}^\XR\fii_j,\om\ra.\end{equation}
In particular, we conclude from \eqref{8} that
\begin{equation}\label{10}
\nabla^\Om_{(v,a)}=\nabla^\Om_v+\nabla_a^\Om.
\end{equation}

The expression of $\nabla_{(v,a)}^{\Omega}$ on smooth cylinder functions
motivates the following definition.

\begin{define}\label{def0}\rom{ The tangent space $T_\om\big(\OM\big)$ to the marked
    configuration space $\OM$ at a point $\om=(\gamma,s)\in\OM$ is defined
    as the Hilbert space
\begin{align*}
T_\om\big(\OM\big):&=L^2(X\to T(X)\dotplus\R;\gamma)\\
&=L^2(X\to T(X);\gamma)\oplus L^2(X\to\R;\gamma)
\end{align*}
with  scalar product
\begin{equation}\label{9.1}
\la
V_\om^1,V_\om^2\ra_{T_\om\left(\OM\right)}
=\int_X\big(\la V_\om^1(x)_{T_x(X)},V_\om^2(x)_{T_x(X)}\ra_{T_x(X)}+
V_\om^1(x)_\R V_\om^2(x)_\R\big)\,\gamma(dx),\end{equation}
where $V_\om^1,V_\om^2\in T_\om\big(\OM\big)$ and $V_\om(x)_{T_x(X)}$ and  $V_\om(x)_\R$
denote the projection of \linebreak $V_\om(x)\in T_x(X)\dotplus\R$ onto $T_x(X)$ and $\R$,
respectively.
The corresponding tangent bundle is
\[ T\big(\OM\big)=\bigcup_{\om\in\OM}T_\om\big(\OM\big).\]
}\end{define}

As usual in Riemannian geometry, having directional derivatives and a Hilbert
space as a tangent space, we can introduce a gradient.

\begin{define}\label{def1}
\rom{
We define an intrinsic gradient of a function $F\colon \OM\to\R$
as a mapping
\[\OM\ni\om\mapsto(\nabla^\Omega F)(\om)\in T_\om\big(\OM\big)\]
such that, for any $(v,a)\in\got$,
\[(\nabla^\Omega_{(v,a)}F)(\om)=\la(\nabla^\Omega F)(\om),(v,a)\ra_{T_\om\left(\OM\right)}.\]
}\end{define}

By \eqref{8} and \eqref{9} we have, for an arbitrary $F\in\FC$ of the form
\eqref{5} and each $\om=(\gamma,s)\in\OM$,
\begin{equation}\label{gradient}(\nabla^\Omega F)(\om;x)=\sum_{j=1}^N\frac{\di g_F}{\di
  r_j}(\la\fii_1,\om\ra,\dots,\la\fii_N,\om\ra)
\nabla^\XR\fii_j(x,s_x),\qquad x\in\gamma.\end{equation}

\begin{rem}\rom{\label{rem3.1}
The operator $s\frac \di{\di s}$ relates to a representation of
the group of dilations
in $L^2(\Rp)$ (cf.\ e.g.\ \cite{Goldin}). Namely, for any $\theta>0$, we put
\[ L^2(\Rp)\ni f(\cdot)\mapsto f(\theta\cdot)\in L^2(\Rp).\]
Setting here $\theta_t^a:=e^{at}$, $a\in \R$, $t\in\R$, we obtain
\begin{equation}\label{13}
\frac d{dt}f(e^{at}s)\big|_{t=0}=sf'(s) a=:(\nabla_a^\Rp f)(s)=(\nabla^\Rp f,a)_{\R}
.\end{equation}

It is easy to verify that the following operator equality holds:
\[\nabla^{\Rp}=S_{\cal L}\nabla^\R S_{\cal L}^{-1},\]
where the mapping
\[ C^1(\R)\ni f\mapsto S_{\cal L}f:=f\circ{\cal L}\in C^1(\Rp)\]
is the isomorphism generated by the bijection $\R_+\ni s\mapsto{\cal L}s=\log
 s\in\R$.
 Notice that, if ${\cal U}_{\cal L}$ is the unitary
between the spaces $L^2(\R)$ and $L^2(\Rp)$ generated by $\cal L$, i.e., if
\[ ({\cal U}_{\cal L}f)(s):=\frac{f({\cal L}s)}{\sqrt{s}}\, ,\qquad f\in
L^2(\R),\, s\in\Rp,\]
then ${\cal U}_{\cal L}\nabla^\R{\cal U}_{\cal L}^{-1}$ coincides with
$\nabla^{\Rp}+\frac12
\pmb{1}$.

Finally, we note that the condition \eqref{trinidad} means, in fact, that an
$\Rp$-derivative of an arbitrary order of a function from $\frak D$ belongs
again to $\frak D$.
}\end{rem}

\subsection{Integration by parts and divergence on the marked Poisson space}

Let the marked configuration space $\OM$ be equipped with the marked Poisson
measure $\pi_{\tsigma}$, and in addition to the condition \eqref{c2} we will
suppose that
\begin{equation}\label{condition}\sqrt q\in H_0^{1,2}(\XR).\end{equation}
Here, $H_0^{1,2}(\XR)$ denotes the local Sobolev space of order 1 constructed with respect
to the gradient $\nabla^\XR$ in the space $L^2_{\mathrm loc}(X;m)\otimes
L^2(\R_+;\lambda)$, i.e., $H_0^{1,2}(\XR)$ consists of functions $f$ defined on
$\XR$ such that, for any set $A\in{\cal B}_{\mathrm c}(\XR)$, the restriction of
$f$ to $A$ coincides with the restriction to $A$ of some function $\fii$ from the
Sobolev space $H^{1,2}(\XR)$ constructed as the closure of $\frak D$ with
respect to the norm
\[ \|\fii\|^2_{1,2}:=\int_{\XR}\Big( |\nabla^X\fii(x,s)|^2_{T_x(X)}+s^2\Big|\frac
{\di}{\di s}\fii(x,s)\Big|^2 +|\fii(x,s)|^2\Big)m(dx)\,\lambda(ds).\]

The set $\FC$ is a dense subset in the space
\[ L^2(\OM,{\cal B}(\OM),\pi_{\tsigma})=:L^2(\pi_{\tsigma}).\]
For any $(v,a)\in\got$, we have a differential operator
 in $L^2(\pi_\tsigma)$ on the domain \linebreak$\FC$ given by
\[\FC\ni F\mapsto\nabla_{(v,a)}^\Omega F\in L^2(\pi_\tsigma).\]
Our aim now is to compute the adjoint operator $\nabla_{(v,a)}^{\Omega\,*}$ in
$L^2(\pi_\tsigma)$. It corresponds, of course, to an integration by parts
formula with respect to the measure $\pi_\tsigma$.

But first we present the corresponding formula on $\XR$.

\begin{define}\rom{
For any $(v,a)\in\got$, the logarithmic derivative of the measure $\tsigma$
along $(v,a)$ is defined as the following function on $\XR$:
\[\beta_{(v,a)}^\tsigma:=\beta_v^\tsigma+\beta_a^\tsigma\]
with
\begin{align*}
\beta_v^\tsigma(x,s):&=\frac{\nabla_v^Xq(x,s)}{q(x,s)}+\div
^Xv(x)\\
&=\bigg\la\frac{\nabla^Xq(x,s)}{q(x,s)},v(x)\bigg\ra_{T_x(X)}+
\div^Xv(x),
\end{align*}
$\div^X=\div^X_m$ being the divergence on $X$ w\rom.r\rom.t\rom.\ $m$\rom, and
\[\beta_a^\tsigma(x,s)
=\frac{s\frac \di{\di s}q(x,s)}{q(x,s)}\,a(x)+a(x)
=\frac{s\frac\di{\di s}p(x,s)}{p(x,s)}\,a(x)+a(x).\]
}\end{define}

Upon \eqref{condition}, we conclude that, for each $(v,a)\in\got$, the function
$\nabla_{(v,a)}^\XR\log q$ is quadratically integrable with respect to the
measure $\tsigma$ on any $\Lext\in{\cal B}_0(\XR)$. Therefore, taking to
notice that $\tsigma(\Lext)<\infty$, we conclude that each function
$\beta^\tsigma_{(v,a)}$, $(v,a)\in\got$, is absolutely integrable on every $\Lext$ with
respect to $\tsigma$.

\begin{lem}[Integration by parts formula on $\XR$]
For all $\fii_1$\rom, $\fii_2\in{\frak D}$\rom, we have
\begin{align*}
&\int_\XR(\nabla^\XR_{(v,a)}\fii_1)(x,s)\fii_2(x,s)\,\tsigma(dx,ds)=\\
&\qquad=-\int_\XR\fii_1(x,s)(\nabla^\XR
_{(v,a)}\fii_2)(x,s)\,\tsigma(dx,ds)
\\&\qquad\quad
\text{}-\int_\XR\fii_1(x,s)\fii_2(x,s)\beta^\tsigma_{(v,a)}(x,s)\,\tsigma(dx,ds).\end{align*}
 \end{lem}

\noindent{\it Proof}. Since $\beta_v^\tsigma$ is exactly the logarithmic
derivative of the measure $\tsigma$ along the vector field $v$, to prove the
lemma it suffices to see that, for arbitrary $f_1,f_2\in C^\infty_{\mathrm
  b}(\Rp)$, $a\in\R$, and for $m$-a.a.\ $x\in X$, we have by virtue of \eqref{c2}
\begin{gather*}
\int_\Rp(\nabla_a^\Rp f_1)(s)f_2(s)\,p(x,ds)=\int_{\Rp}asf_1'(s)f_2(s)
p(x,s)\,\lambda(ds)\\
=af_1(s)f_2(s)sp(x,s)\Big|_{s=0}^{s=\infty}-\int_\Rp f_1(s)a\,\frac
d{ds}\big(f_2(s)sp(x,s)\big)\,\lambda(ds)\\
=-\int_{\Rp}f_1(s)(\nabla_a^\Rp f_2)(s)\,
p(x,ds)-\int_\Rp f_1(s) f_2(s)\beta_a^{p(x,\cdot)}(s)\,p(x,ds).\end{gather*}
Here
\[\beta_a^{p(x,\pmb{\cdot})}(s)=\frac{s\frac \di{\di s}p(x,s)}{p(x,s)}\,a+a\]
and  we have used the formulas
\eqref{13} and \eqref{c2}.\quad $\blacksquare$

\begin{define}\rom{
For any $(v,a)\in\got$, the logarithmic derivative of the marked Poisson measure
$\pi_\tsigma$ along $(v,a)$ is defined as the following function on $\OM$:
\begin{equation}\label{16}\OM\ni\om\mapsto
B^{\pi_\tsigma}_{(v,a)}(\om):=\la\beta_{(v,a)}^\tsigma,\om\ra.\end{equation}
}\end{define}

A motivation for this definition is given by the following theorem.

\begin{th}
[Integration by parts formula]
For all $F,G\in\FC$ and each $(v,a)\in\got$\rom, we have
\begin{align}
\int_{\OM}(\nabla_{(v,a)}^\Omega F)(\om)G(\om)\,\pi_\tsigma(d\om)&=-\int_{\OM}
F(\om)(\nabla_{(v,a)}^\Omega G)(\om)\,\pi_\tsigma(d\om)\notag\\
&\quad-\int_{\OM} F(\om) G(\om) B_{(v,a)}^{\pi_\tsigma}(\om)\,\pi_\tsigma(d\om),
\label{11}\end{align}
or
\begin{equation}
\label{12}
\nabla_{(v,a)}^{\Om\,*}=-\nabla_{(v,a)}^\Omega-B_{(v,a)}^{\pi_{\tsigma}}(\om)\end{equation}
as an operator equality on the domain $\FC$ in $L^2(\pi_\tsigma)$\rom.
\label{thIbP}\end{th}

\noindent{\it Proof}. Because of \eqref{10}, the formula \eqref{12} will be
proved if we prove it first for the operator $\nabla^\Om_v$, i.e., when $a(x)\equiv
0$, and then for the operator $\nabla^\Om_a$, i.e., when $v(x)=0\in T_x(X)$ for all
$x\in X$. We present below only the proof for $\nabla^\Om_a$, since the proof for
$\nabla^\Om_v$ is basically the same as that of the integration by parts
formula in case of Poisson measures \cite{AKR}.


By Proposition~\ref{prop2}, we have for all $a\in C^\infty_0(X)$
\[\int_{\OM} F(\theta_t^a(\om))G(\om)\,\pi_\tsigma(d\om)=\int_{\OM}
F(\om)G(\theta_{-t}^a(\om))\,
\pi_{\theta_t^{a\,*}\tsigma}(d\om).\]
Differentiating this equation with respect to $t$, interchanging $d/dt$ with the
integrals and setting $t=0$, the l.h.s.\ becomes the l.h.s.\ of \eqref{11}. To see
that the r.h.s.\ then also coincides with the r.h.s.\ of \eqref{11}, we note
that
\[\frac d{dt}G(\theta_{-t}^a(\om))\big|_{t=0}=-(\nabla_a^\Omega G)(\om),\]
and by Proposition~\ref{prop3}
\begin{gather*}
\frac
d{dt}\bigg[\frac{d\pi_{\theta_t^{a\,*}\tsigma}}{d\pi_\tsigma}(\om)\bigg]\bigg|_{t=0}=
\sum_{(x,s)\in\om}\frac d{dt}\, p_{\theta_t^a}^\tsigma (x,s)
\bigg|_{t=0}\\
=\sum_{(x,s)\in\om}\frac
d{dt}\bigg[\frac{p(x,e^{-ta(x)}s)}{p(x,s)e^{ta(x)}}\bigg]\bigg|_{t=0}
=-\la\beta_a^\tsigma,\om\ra
=-B_a^{\pi_\tsigma}(\om).\quad\blacksquare
\end{gather*}

\begin{define}\label{def2}\rom{
For a vector field
\[V\colon\OM\ni\om\mapsto V_\om\in T_\om(\OM),\]
the divergence $\div_{\pi_\tsigma}^\Omega V$ is defined via the duality relation
\[ \int_{\OM}\la V_\om,\nabla^\Omega
F(\om)\ra_{T_\om\big(\OM\big)}\pi_\tsigma(d\om)=-\int_{\OM}F(\om)
(\div_{\pi_\tsigma}^\Omega V)(\om)\,\pi_\tsigma(d\om)\]
for all $F\in\FC$, provided it exists (i.e., provided
\[ F\mapsto\int_{\OM}\la V_\om,\nabla^\Omega
F(\om)\ra_{T_x\big(\OM\big)}\pi_\tsigma(d\om)\]
is continuous on $L^2(\pi_\tsigma)$).
}\end{define}

A class of smooth vector fields on $\OM$ for which the divergence can be
computed in an explicit form is described in the following proposition.

\begin{prop}\label{prop4} For any vector field
\[ V_\om(x)=\sum_{j=1}^N G_j(\om)(v_j(x),a_j(x)),\qquad \om\in\OM,\ x\in X,\]
with $G_j\in\FC$\rom, $(v_j,a_j)\in\got$\rom, $j=1,\dots,N$\rom, we have
\begin{align*}
(\div^\Omega_{\pi_\tsigma}V)(\om)&=\sum_{j=1}^N\big(\nabla_{(v_j,a_j)}^\Omega
G_j\big)(\om)+\sum_{j=1}^N B_{(v_j,a_j)}^{\pi_{\tsigma}}(\om)G_j(\om)\\
&=\sum_{j=1}^N
\la\nabla^\Omega G_j(\om), (v_j,a_j)\ra_{T_\om\big(\OM\big)}+\sum_{j=1}^N\la\beta^\tsigma_{(v_j,a_j)}
,\om\ra G_j(\om).\end{align*}\end{prop}

\noindent{\it Proof}. Due to the linearity of $\nabla^\Omega$, it is sufficient
to consider the case $N=1$, i.e., $V_\om(x)=G(\om)(v(x),a(x))$. By
Theorem~\ref{thIbP}, we have for all $F\in\FC$
\begin{gather*}
-\int_{\OM}\la V_\om,\nabla^\Omega F(\om)\ra_{T_\om\big(\OM\big)}\pi_{\tsigma}(d\om)=
-\int_{\OM}G(\om)\nabla_{(v,a)}^\Omega F(\om)\,\pi_\tsigma(d\om)\\
=\int_{\OM}\big(\nabla^\Omega_{(v,a)}G\big)(\om)F(\om)\,\pi_\tsigma(d\om)+\int
_{\OM}G(\om)F(\om)B_{(v,a)}^{\pi_\tsigma}(\om)\,\pi_{\tsigma}(d\om),\end{gather*}
which yields
\begin{align*}
(\div_{\pi_\tsigma}^\Omega
V)(\om)&=\nabla^\Omega_{(v,a)}G(\om)+B^{\pi_\tsigma}_{(v,a)}(\om)G(\om)\\
&=\la\nabla^\Omega G(\om),(v,a)\ra_{T_\om\big(\OM\big)}+\la\beta^\tsigma_{(v,a)},\om\ra
G(\om).\quad \blacksquare\end{align*}

\begin{rem}\label{rem1}\rom{
Extending the definition of $B^{\pi_\tsigma}$ in \eqref{16} to the class of
vector fields $V=\sum_{j=1}^N G_j\otimes (v_j,a_j)$ by
\[ B_V^{\pi_\tsigma}(\om):=\sum_{j=1}^N\la\beta^\tsigma_{(v_j,a_j)},\om\ra
G_j(\om)+\sum_{j=1}^N
\big(\nabla^\Omega_{(v_j,a_j)}G_j\big)(\om),\]
we obtain that
\[\div^\Omega_{\pi_\sigma}\cdot=B^{\pi_\tsigma}_{\pmb{\pmb{.}}}.\]
In particular, if $(v,a)\in\got$, it follows,  for the ``constant'' vector field
$V_\om\equiv(v,a)$ on $\OM$, that
\[ \div^\Omega_{\pi_\tsigma}(v,a)(\om)=\la\div^{\XR}_\tsigma(v,a),\om\ra,\]
where $\div^\XR_\tsigma(v,a)=\beta^\tsigma_{(v,a)}$ is the divergence of
$\tsigma$ on $\XR$ w.r.t.\ $(v,a)$:
\begin{align*}
&\int_{\XR}\la\nabla^\XR\fii(x,s),(v(x),a(x))\ra_{T_{(x,s)}(\XR)}\,\tsigma(dx,ds)\\
&\qquad=-\int_\XR\fii(x,s)\big(\div_\tsigma^{\XR}(v,a)\big)(x,s)\,\tsigma(dx,ds),\qquad
\fii\in{\frak D}.\end{align*}
}\end{rem}

\subsection{Integration by parts characterization}

In the works \cite{AKR,AKR2} it was shown that the mixed Poisson measures are
exactly the ``volume elements'' corresponding to the differential geometry on the
configuration space $\Gamma_X$. Now, we wish to prove that an analogous
statement holds true in our case of $\OM$ for mixed marked Poisson measures.

We start with a lemma that describes $\tsigma$ as a unique (up to a constant)
measure on $\XR$ with respect to which the divergence $\div^\XR_\tsigma$ is the
dual operator of the gradient $\nabla^{\XR}$.

\begin{lem}\label{lemcharacter}
Let
\begin{align}
\frac{\nabla^X_vq(x,s)}{q(x,s)}&\in L^1_{\mathrm
  loc}(X\times\Rp;m\otimes\lambda), \qquad v\in V_0(X),\notag\\
\frac{\frac{\di}{\di s}p(x,s)}{p(x,s)}&\in L^1_{\mathrm loc}(X\times \Rp;
  m\otimes\lambda).\label{local}
\end{align}
Then\rom, for every $\Lambda\in\Oc(X)$ the measures $z\tsigma$\rom, $z\ge
0$\rom, are the only positive Radon measures $\xi$ on $\Lext$ such that
$\div_\tsigma^\XR$ is the dual operator on $L^2(\xi)$ of $\nabla^\XR$ when
considered with domains $V_0(\Lambda)\times C^\infty_0(\Lambda)$\rom, resp\rom.\
$C^\infty_{\text{\rom{0,b}}}(\Lext)$
\rom(i\rom.e\rom., the set of all $(v,a)\in\got$\rom, resp\rom.\
$\varphi\in{\frak D}$ with support in $\Lambda$\rom, resp\rom.\ $\Lext$\rom{).}
\end{lem}

\noindent {\it Proof}. By using the condition \eqref{local},
the lemma is
obtained in complete analogy with Remark~4.1 (iii) in \cite{AKR2}. Indeed, let
$q_1(x,s)$ and $q_2(x,s)$ be two densities w.r.t.\ $m\otimes \lambda$ for which
the logarithmic derivatives coincide. Then, we get
\begin{align*}
\nabla_v^X\log q_1(x,s)&=\nabla_v^X\log q_2(x,s),\qquad v\in V_0(X),\\
\frac\di {\di s}\log q_1(x,s)&=\frac{\di }{\di s}\log q_2(x,s)\quad
\text{$m\otimes\lambda$-a.s.},
\end{align*}
which yields respectively
\begin{align*}q_1(x,s)&=q_2(x,s)c(s),\\
q_1(x,s)&=q_2(x,s)\tilde c(x)\quad\text{$m\otimes\lambda$-a.s.}
\end{align*}
Therefore, $q_1(x,s)=\operatorname{const} q_2(x,s)$ $m\otimes\lambda$-a.s.\quad $\blacksquare$

Let $\nu$ be a probability measure on $(\R_+,{\cal B}(\Rp))$. Then, we define a
mixed Poisson measure as follows:
\begin{equation}\label{salko} \mu_{\nu,\tsigma}=\int_\Rp\pi_{z\tsigma}\,\nu(dz).\end{equation}
Here, $\pi_{0\tsigma}$ denotes the Dirac measure on $\OM$ with mass in $\om=\varnothing$.
 Let ${\cal M}_l(\OM)$, $l\in[1,\infty)$, denote the set of all probability
 measures on $(\OM,{\cal B}(\OM))$ such that
\[\int_{\OM}|\la f,\om\ra|^l\,\mu(d\om)<\infty\quad \text{for all }f\in
 C_{\text{0,b}}(\XR),\ f\ge0.\]
Clearly, $\mu_{\nu,\tsigma}\in {\cal M}_l(\OM)$ if and only if
\begin{equation}\label{30}
\int_{\Rp}z^l\,\nu(dz)<\infty.\end{equation}
We define $(\text{IbP})^\tsigma$ to be the set of all
$\mu\in {\cal
  M}_1(\OM)$ with the property that $\om\mapsto \la\beta^\tsigma_{(v,a)},\om\ra$
is $\mu$-integrable for all $(v,a)\in\got$ and which satisfy \eqref{11} with $\mu$
replacing $\pi_\tsigma$ for all $F,G\in\FC$, $(v,a)\in\got$.
 We note that \eqref{11} makes sense only for such measures and that
 $B_{(v,a)}^{\pi_\tsigma}$ depends only on $\tsigma$ not on
 $\pi_\tsigma$. Obviously, since $\nabla^\XR_{(v,a)}$ obeys the product rule
 for all $(v,a)\in\got$, we can always take $G\equiv
 1$. Furthermore, $(\text{IbP})^{\tsigma}$ is convex.

\begin{th}\label{hahaha} Let the condition \eqref{local} be satisfied\rom. Then\rom, the
  following conditions are equivalent\rom:\\
\rom{(i)} $\mu\in({\mathrm IbP})^\tsigma$\rom;\\
\rom{(ii)} $\mu=\mu_{\nu,\tsigma}$ for some probability measure $\nu$ on
$(\Rp,{\cal B}(\Rp))$ satisfying \eqref{30} with $l=1$\rom. \end{th}

\noindent{\it Proof}. The part (ii)$\Rightarrow$(i) is trivial. In order to
prove (i)$\Rightarrow$(ii),
we have to modify the corresponding part of the proof of Theorem~4.1
in \cite{AKR2}. We present this modification in Appendix.\vspace{2mm}\quad $\blacksquare$

As a direct consequence of Theorem~\ref{hahaha}, we obtain

\begin{cor} Suppose that the condition \eqref{local} is
  satisfied\rom. Then the extreme points of $(\operatorname{IbP})^\tsigma$ are
  exactly $\pi_{z\tsigma}$\rom, $z\ge0.$
\end{cor}

\subsection{A lifting of the geometry}

Just as in the case of the geometry on the configuration space, we can present
an interpretation of the formulas obtained in subsections~3.1--3.3 via a simple
``lifting rule.''

Suppose that $f\in C_{\text{0,b}}(\XR)$, or more generally $f$ is an arbitrary  measurable
function on $\XR$ for which there exists (depending on $f$) $\Lambda\in\Bc(X)$
such that $\supp f\subset\Lext$. Then, $f$  generates a (cylinder) function on $\OM$
by the formula
\[ L_f(\om):=\la f,\om\ra,\qquad \om\in\OM.\]
We will call $L_f$ the lifting of $f$.

As before, any vector field $(v,a)\in\got$,
\[(v,a)\colon X\ni x\mapsto (v(x),a(x))\in T_{(x,s)}(\XR)=T_x(X)\dotplus\R,\]
can be considered as a vector field on $\OM$ (the lifting of $(v,a)$), which we
denote by $L_{(v,a)}$:
\[L_{(v,a)}\colon \OM\ni\om=\{\gamma,s\}\mapsto\{x\mapsto(v(x),a(x))\}\in
T_\om(\OM)=L^2(X\to T(X)\dotplus\R\,;\gamma).\]
For $(v,a), (w,b)\in\got$, the formula \eqref{9.1} can be written as follows:
\[\big\la L_{(v,a)},L_{(w,b)}\big\ra_{T_\om\big(\OM\big)}=L_{\la(v,a),(w,b)\ra_{T(\XR)}}(\om),\]
i.e., the scalar product of lifted vector fields is computed as the lifting of
the scalar product
\[\la(v(x),a(x)),(w(x),b(x))\ra_{T_{(x,s)}(\XR)}=f(x).\]
This rule can be used as a definition of the tangent space $T_\om\big(\OM\big)$.

The formula \eqref{8} has now the following interpretation:
\begin{equation}\label{20}
\big(\nabla^\Omega_{(v,a)}L_\fii\big)(\om)=L_{\nabla^\XR_{(v,a)}\fii}(\om),\qquad
\fii\in{\frak D},\ \om\in\OM,
\end{equation}
and the ``lifting rule'' for the gradient is given by
\begin{equation}\label{17}(\nabla^\Omega L_\fii)(\gamma,s)\colon\gamma\ni
  x\mapsto\nabla^{\XR}\fii( x,s_x).\end{equation}

As follows from \eqref{16}, the logarithmic derivative $B^{\pi_\tsigma}_{(v,a)}\colon\OM\to\R$
is obtained via the lifting procedure of the corresponding logarithmic
derivative $\beta^\tsigma_{(v,a)}\colon\XR\to\R$, namely,
\[ B^{\pi_\tsigma}_{(v,a)}(\om)=L_{\beta^\tsigma_{(v,a)}}(\om),\]
or equivalently, one has for the divergence of a lifted vector field:
\begin{equation}\label{18}\div^\Omega_{\pi_\tsigma}L_{(v,a)}=L_{\div^\XR_\tsigma}(v,a).\end{equation}

We underline that by \eqref{20} and \eqref{17} one recovers the action of
$\nabla^\Omega_{(v,a)}$ and $\nabla^\Omega$ on all functions from $\FC$
algebraically from requiring the product or the chain rule to hold. Also, the
action of $\div_{\pi_\tsigma}^\Omega$ on more general cylindrical vector fields
follows as in Remark~\ref{rem1} if one assumes the usual product rule for
$\div_{\pi_\tsigma^\Omega}$ to hold.

\subsection{Representations of the Lie algebra of the group $\Got$}

Using the $\Got$-quasiinvariance of $\pi_{\tsigma}$, we can define a unitary
representation of the group $\Got=\Diff\underset{\alpha}{\times}\R_+^X$ in the
space $L^2(\pi_{\tsigma})$. Namely, for $g\in\Got$, we define a unitary operator
\[\big(V_{\pi_\tsigma}(g)F\big)(\om):=F(g(\om))\sqrt{\frac{dg^{-1*}\pi_\tsigma
    }{d\pi_\tsigma}(\om)},\qquad F\in L^2(\pi_\tsigma).\]
Then, we have
\[ V_{\pi_\tsigma}(g_1)V_{\pi_\tsigma}(g_2)=V_{\pi_\tsigma}(g_1g_2),\qquad
g_1,g_2\in\Got.\]

As has been noted in the introduction,
we have the following proposition.

\begin{prop}\label{propreduce}
The representation $V_{\pi_\tsigma}$ of
$\Got$ is reducible\rom.
\end{prop}

\noindent{\it Proof}.
Let us suppose for simplicity of notations that $\tsigma$ is a
product measure: $\tsigma=\hat\sigma=\sigma\otimes\tau$,
$\tau(ds)=p(s)\,\lambda(ds)$ (it will be seen that our proof can be generalized
to the case of a general $\tsigma$ without difficulties).
Consider the regular representation of the
group of dilations on $L^2(\tau)$:
\[ (U_\tau(\theta)\fii)(s)=\fii(\theta s)\sqrt{\frac{p(\theta s)}{p(s)}},\qquad
\theta\in\Rp,\, \fii\in L^2(\tau).\]
This representation is reducible (one can also easily
see that it is unitarily equivalent to the regular
representation of the additive group $\R$ on $L^2(\R)$, see e.g.\
\cite{VilKlimyk}). Let $\cal H$ be an arbitrary subspace of $L^2(\tau)$ that is
invariant with respect to this representation and such that  ${\cal H}\ne\{0\}$
and  ${\cal H}\ne
L^2(\tau)$. Then, the subspace $L^2(\sigma)\otimes{\cal H}$ of $L^2(\tsigma)$ is
invariant with respect to the regular representation of the group $\Got $ on
$L^2(\tsigma)$:
\[(V_{\tsigma}(g)f)(x,s)=f(g(x,s))\sqrt{\frac{dg^{-1*}\tsigma}{d\tsigma}(x,s)},\qquad
g\in\Got,\, f\in L^2(\tsigma).\]

For a set $\Upsilon\in{\cal B}(\XR)$ denote by ${\cal B}_\Upsilon(\OM)$ the
$\sigma$-algebra on $\OM$ generated by the mappings $N_B$ (defined by
\eqref{NB}),
where $B\in\Bc(\XR)$
and $B\subset\Upsilon$. By definition, each $B_\Upsilon(\OM)$ is a
sub-$\sigma$-algebra of ${\cal B}(\OM)$. It is easy to see that, for any
$\Lambda\in\Oc(X)$, there is a natural isomorphism between the $\sigma$-algebras
${\cal B}_{\Lext}(\OM)$ and ${\cal B}(\Om_\Lambda^\Rp)$. Then, by the definition of a marked
Poisson measure, we get
\begin{gather*}
L^2_{\Lambda}(\pi_\tsigma):=L^2(\OM,{\cal B}_{\Lext}(\OM),\pi_\tsigma)\backsimeq
L^2(\Om_\Lambda^\Rp,{\cal B}(\Om_\Lambda^\Rp),\pi^\Lambda_\tsigma)\\
=\bigoplus_{n=0}^\infty L^2(\Om^{\Rp}_{\Lambda}(n),{\cal
  B}(\Om_\Lambda^\Rp(n)),\pi_\tsigma^\Lambda \restriction \Om_\Lambda^\Rp(n))\\
=e^{-\tsigma(\Lext)}\bigoplus_{n=0}^\infty \frac 1{n!}\,
L^2(\Om_\Lambda^\Rp(n),{\cal B}(\Om_\Lambda^{\Rp}(n)),\tsigma_{\Lambda,n})\\
\backsimeq e^{-\tsigma(\Lext)}\bigoplus_{n=0}^\infty \frac 1{n!}\,\hat
L^2(\Lambda_{\mathrm\ext}^n,\tsigma^{\otimes n}),\end{gather*}
where $\hat L^2(\Lambda_{\mathrm \ext}^n,\tsigma^{\otimes n})=(L^2(\Lext,\tsigma))^{\widehat \otimes
  n}$, $\widehat \otimes$ denoting symmetric tensor product.

We note that $L^2_{\Lambda}(\pi_\tsigma)$ is a subspace of
$L^2_{\Lambda'}(\pi_\tsigma)$ provided $\Lambda\subset\Lambda'$, and the
whole space $L^2(\pi_\tsigma)$ is the ``limit'' of the spaces
$L^2_{\Lambda}(\pi_\tsigma)$ as $\Lambda\nearrow X$, i.e., each function $F\in
L^2(\pi_\tsigma)$ can be represented as the $L^2(\pi_\tsigma)$ limit of a
sequence of functions $F_\Lambda\in L^2_{\Lambda}(\pi_\tsigma)$ as $\Lambda \nearrow
X$.

Now, for each $\Lambda\in\Oc(X)$, we define a subspace $R_\Lambda$ of
$L^2_{\Lambda}(\pi_\tsigma)$ as follows:
\begin{align*}
R_\Lambda:&=e^{-\tsigma(\Lext)}\bigoplus_{n=0}^\infty\frac1{n!}\,(L^2(\Lambda;\sigma)\otimes
  {\cal H })^{\widehat\otimes n}\\
&=e^{-\tsigma(\Lext)}\bigoplus_{n=0}^\infty\frac 1{n!}\,\hat
  L^2(\Lambda^n;\sigma^{\otimes n})\otimes{\cal H}^{\widehat \otimes n}.\end{align*}
Again, for arbitrary $\Lambda,\Lambda'\in\Oc(X)$, $\Lambda\subset \Lambda'$, we
  obtain the inclusion $R_\Lambda\subset R_{\Lambda'}$, and let $R$ be the
  ``limit'' of the $R_\Lambda$ spaces as $\Lambda \nearrow X$. Evidently, $R$ does
  not coincide with the whole $L^2(\pi_\tsigma)$ (an arbitrary function from
  the space
\[e^{-\tsigma(\Lext)}\bigoplus_{n=0}^\infty \frac1{n!}\,
 \hat  L^2(\Lambda^n;\sigma^{\otimes n})\otimes({\cal H} ^\bot )^{\widehat\otimes n},\]
where ${\cal H} ^\bot$ is the orthogonal complement of ${\cal H}$ in $L^2(\tau)$,
will be orthogonal to all space $R_{\Lambda'}$ with $\Lambda\subset\Lambda'$,
and therefore to $R$).

It is easy to see that, for each fixed $g\in\Got$, all spaces $R_\Lambda$ such
that $K_g\subset\Lambda$ are invariant with respect to the operator
$V_{\pi_\tsigma}(g)$, and hence so is the space $R$, which concludes the
proof.\quad\vspace{2mm} $\blacksquare$

As in subsec.~3.1, to any vector field $v\in V_0(X)$ there corresponds a
one-parameter subgroup of diffeomorphisms $\psi_t^v$, $t\in\R$. It generates a
one-parameter unitary group
\[V_{\pi_\tsigma}(\psi_t^v):=\exp[itJ_{\pi_\tsigma}(v)],\qquad t\in\R,\]
where $J_{\pi_\tsigma}(v)$ denotes the self adjoint generator of this
group. Analogously, to  a one-parameter unitary group $\theta_t^a$, $a\in
C^\infty_0(X)$, there corresponds a one-parameter unitary group
\[ V_{\pi_\tsigma}(\theta_t^a):=\exp [it I _ {\pi_\tsigma}(a)]\]
with a generator $I_{\pi_\tsigma}(a)$.

\begin{prop}\label{prop5}
For any $v\in V_0(X)$ and $a\in C^\infty_0(X)$\rom, the following operator equalities
on the domain $\FC$ hold\rom:
\begin{align*}
J_{\pi_\tsigma}(v)&=\frac 1i\,\nabla_v^\Omega+\frac{1}{2i}\,B_v^{\pi_\tsigma},\\
I_{\pi_\tsigma}(a)&=\frac1i\,\nabla_a^\Omega+\frac1{2i}\,B_a^{\pi_\tsigma}.
\end{align*}\end{prop}

\noindent{\it Proof}. These equalities follow immediately from the definition of
the directional derivatives $\nabla_v^\Omega$ and $\nabla_a^\Omega$,
Theorem~\ref{thIbP}, and the form of the operators $V_{\pi_\tsigma}(\psi_t^v)$
and $V_{\pi_\tsigma}(\theta_t^a)$.\quad $\blacksquare$ \vspace{2mm}

For any $(v,a)\in\got$, define an operator
\[{\cal R}_{\pi_\tsigma}(v,a):= J_\pit(v)+I_\pit(a).\]
By Proposition~\ref{prop5},
\[{\cal R}_\pit(v,a)=\frac1i\nabla_{(v,a)}^\Om+\frac1{2i}B^\pit_{(v,a)}.\]
We wish to derive now a commutation relation between these operators.

\begin{lem}\label{lem3.3} The Lie-bracket $[(v_1,a_1),(v_2,a_2)]$ of the vector fields
  $(v_1,a_1)$\rom, $(v_2,a_2)\in \got$\rom, i\rom.e\rom{.,} a
  vector field from $\got$ such that
\[
  \nabla^{\XR}_{[(v_1,a_1),(v_2,a_2)]}=\nabla^{\XR}_{(v_1,a_1)}\nabla^\XR_{(v_2,a_2)}-
\nabla^\XR_{(v_2,a_2)}\nabla^\XR_{(v_1,a_1)}\quad \text{\rom{on} }{\frak D},\]
is given by
\[[(v_1,a_1),(v_2,a_2)]=([v_1,v_2],\nabla^X_{v_1}a_2-\nabla^X_{v_2}a_1),\]
where $[v_1,v_2]$ is the Lie-bracket of the vector fields $v_1$\rom, $v_2$ on
$X$\rom.
\end{lem}

\noindent{\it Proof}. The lemma is obtained by a direct computation if one uses
the following evident relations which hold on ${\frak D}$:
\begin{alignat*}{2}
\nabla^X_{v_1}\nabla^X_{v_2}-\nabla^X_{v_2}\nabla^X_{v_1}&=\nabla^X_{[v_1,v_2]},&&
\qquad v_1, v_2\in
V_0(X),\\
\nabla^{\Rp}_{a_1}\nabla_{a_2}^{\Rp}-\nabla_{a_2}^{\Rp}\nabla^{\Rp}_{a_1}&=0,&&
\qquad a_1,a_2\in C^\infty_0(X),\\
\nabla^X_{v}\nabla^{\Rp}_a-\nabla_a^{\Rp}\nabla^X_v&=\nabla^{\Rp}_{\nabla^X_v
  a},&&\qquad v\in V_0(X),\ a\in C_0^\infty(X).\quad\blacksquare\end{alignat*}

\begin{prop}\label{prop3.3}
For arbitrary $(v_1,a_1),\, (v_2,a_2)\in\got$\rom, the
following operator equality holds on $\FC$\rom:
\begin{align*}
[{\cal R}_\pit(v_1,a_1),{\cal R}_\pit(v_2,a_2)]&={\cal
  R}_\pit([(v_1,a_1),(v_2,a_2)])\\
&={\cal R}_\pit([v_1,v_2],\nabla^X_{v_1}a_2-\nabla^X_{v_2}a_1).\end{align*}
In particular\rom,
\begin{alignat*}{2}
[J_\pit(v_1),J_\pit(v_2)]&=-iJ_\pit([v_1,v_2]),&&\qquad v_1,v_2\in V_0(X),\\
[I_\pit(a_1),I_\pit(a_2)]&=0,&&\qquad a_1,a_2\in C_0^\infty(X),\\
[J_\pit(v),I_\pit(a)]&=-iI_\pit(\nabla_v^Xa),&&\qquad v\in V_0(X),\, a\in
C_0^\infty(X).\end{alignat*}
\end{prop}

\noindent {\it Proof}. First we note that Lemma~\ref{lem3.3} and \eqref{8}
immediately imply
\[ \nabla^\Omega_{(v_1,a_1)}\nabla^\Omega_{(v_2,a_2)}-\nabla^\Omega_{(v_2,a_2)}\nabla^\Om_{(v_1,a_1)}
=\nabla^\Om_{[(v_1,a_1),(v_2,a_2)]}\quad \text{on }\FC.\]
Therefore, by using the chain rule, we conclude that the lemma will be proved if
we show that
\begin{equation}\label{zyrk}
\nabla^\Om_{(v_1,a_1)}B^\pit_{(v_2,a_2)}-\nabla^\Om_{(v_2,a_2)}B^\pit_{(v_1,a_1)}=B^\pit_{[(v_1,a_1),(v_2,
  a_2)]}\quad \text{$\pit$-a.e.}\end{equation}
But upon the representation
\[ B^\pit_{(v,a)}(\om)=\la\nabla^\XR_{(v,a)}\log q+\div^Xv+a,\om\ra,\]
we easily derive \eqref{zyrk} again from Lemma~\ref{lem3.3}.\quad $\blacksquare$

Thus, the operators ${\cal R}_{\pi_\tsigma}(v,a)$, $(v,a)\in\got$, give a marked
Poisson space representation of the Lie algebra $\got$ of the group $\Got$.

\section{Intrinsic Dirichlet forms\\ on marked Poisson spaces}

\subsection{Definition of the intrinsic Dirichlet form}

We start with introducing some useful spaces of smooth cylinder functions on
$\OM$ in addition to $\FC$. By $\FP$ we denote the set of all cylinder functions
of the form \eqref{5} in which the generating function $g_F$ is a polynomial on
$\R^N$, i.e., $g_F\in{\cal P}(\R^N)$. Analogously, we define $\FCp$ where
$f_F\in C_{\mathrm p}^\infty(\R^N)$ (:=the set of all $C^\infty$-functions $f$
on $\R^N$ such that $f$ and its partial derivatives of any order are
polynomially bounded).

We have obviously
\begin{align*}
\FC&\subset\FCp,\\
\FP&\subset\FCp,\end{align*}
and these are algebras with respect to the usual operations. The existence of
the Laplace transform $\ell_{\pi_\tsigma}(f)$ for each $f\in C_{\mathrm
  0,b}(\XR)$ implies, in particular, that $\FCp\subset L^2(\pi_\tsigma)$.

\begin{define}\rom{
For $F,G\in\FCp$, we introduce a pre-Dirichlet form as
\begin{equation}\label{4.1}
\EOM(F,G)=\int_{\OM}\la\nabla^\Om F(\om),\nabla^\Om
G(\om)\ra_{T_\om(\OM)}\,\pi_{\tsigma}(d\om).\end{equation}
}\end{define}

Note that, for all $F\in\FCp$, the formula \eqref{gradient} is still valid and
therefore, for $F=g_F(\la\fii_1,\cdot\ra,\dots,\la\fii_N,\cdot\ra)$ and
$G=g_G(\la\xi_1,\cdot\ra,\dots,\la\xi_M, \cdot\ra)$ from $\FCp$, we have
\begin{gather}
\la\nabla^\Om F(\om),\nabla^\Om G(\om)\ra_{T_\om(\OM)}=\notag\\
=\sum_{j=1}^N\sum_{k=1}^M\frac{\di g_F}{\di
  r_j}(\la\fii_1,\om\ra,\dots,\la\fii_N,\om\ra)\frac{\di g_G}{\di
  r_k}(\la\xi_1,\om\ra,\dots,\la\xi_M,\om\ra)\times\notag\\
\times\int_X\la\nabla^{\XR}
\fii_j(x,s_x),\nabla^{\XR}\xi_k(x,s_x)\ra_{T_{(x,s_x)}(\XR)}\,\gamma(dx)\notag\\
=\sum_{j=1}^N\sum_{k=1}^M\frac{\di g_F}{\di
  r_j}(\la\fii_1,\om\ra,\dots,\la\fii_N,\om\ra)\frac{\di g_G}{\di
  r_k}(\la\xi_1,\om\ra,\dots,\la\xi_M,\om\ra)\times\notag\\
\times\la\la\nabla^{\XR}\fii_j,\nabla^{\XR}\xi_k\ra_{T(\XR)},\om\ra.\label{sex}
\end{gather}
Since for $\fii,\xi\in{\frak D}$, the function
\begin{gather*}\la\nabla^\XR\fii(x,s),\nabla^\XR\xi(x,s)\ra_{T_{(x,s)}(\XR)}=\\
=\la\nabla^X\fii(x,s),\nabla^X\xi(x,s)\ra_{T_x(X)}+\nabla^{\Rp}\fii(x,s)\nabla^\Rp\xi(x,s)\end{gather*}
belongs to $C_{\mathrm 0,b}(\XR)$ (see \eqref{trinidad}), we conclude that
\[\la\nabla^\Om F(\cdot),\nabla^\Om(\cdot)\ra_{T(\OM)}\in
L^1(\pi_\tsigma),\qquad F,G\in\FCp,\]
and so \eqref{4.1} is well defined.

We will call $\EOM$ the intrinsic pre-Dirichlet form corresponding to the marked
Poisson measure $\pi_\tsigma$ on $\OM$. In the next subsection we will prove the
closability of $\EOM$.

\subsection{Intrinsic Dirichlet operators}

Let us introduce a differential operator $\HOM$ on the domain $\FC$ which is
given on any $F\in\FC$ of the form \eqref{5} by the formula
\begin{gather}
(\HOM F)(\om):=\notag\\
=-\sum_{j,k=1}^N \frac{\di^2 g_F}{\di r_j\di
  r_k}(\la\fii_1,\om\ra,\dots,\la\fii_N,\om\ra)\la\la\nabla^{\XR
  }\fii_j,\nabla^\XR\fii_k\ra_{T(\XR)},\om\ra\notag\\
\mbox{}-\sum_{j=1}^N\frac{\di g_F}{\di
  r_j}(\la\fii_1,\om\ra,\dots,\la\fii_N,\om\ra)\times\notag\\
\times \int_{\XR}\Big[\Delta^X\fii_j(x,s)+s^2\frac{\di^2}{\di
  s^2}\fii_j(x,s)\notag\\
\mbox{}+\la\nabla^{\XR}\log
q(x,s),\nabla^\XR\fii_j(x,s)\ra_{T_{(x,s)}(\XR)}+2s\frac \di{\di
  s}\fii_j(x,s)\Big]\om(dx,ds),
\label{4.3}\end{gather}
where $\Delta^X$ denotes the Laplace--Beltrami operator corresponding to
$\nabla^X$. Since
\[ \la\nabla^\XR\log q,\nabla^\XR\fii_j\ra_{T(\XR)}\in L^2(\tsigma)\cap
L^1(\tsigma)\]
(see subsec.~3.2, in particular, the condition \eqref{condition}), the r.h.s.\
of \eqref{4.3} is well defined as an element of $L^2(\pi_\tsigma)$. To show that
the operator $\HOM$ is well defined, we still have to show that its definition
does not depend on the representation of $F$ in \eqref{5}, which will be done
below.

Let us consider also the pre-Dirichlet operator corresponding to the measure
$\tsigma$ on $\XR$ and to the gradient $\nabla^{\XR}$:
\begin{gather}
{\cal
  E}^\XR_\tsigma(\fii,\xi):=\int_\XR\la\nabla
^\XR\fii(x,s),\nabla^\XR\xi(x,s)\ra_{T_{(x,s)}(\XR)}\,\tsigma (dx, ds)\notag \\
=\int_{\XR}\big[\la\nabla^X\fii(x,s),\nabla^X\xi(x,s)\ra_{T_x(X)}+\nabla^{\Rp}\fii(x,s)\nabla^\Rp\xi(x,s)\big]
  \,\tsigma(dx,ds),\label{4.5}\end{gather}
where $\fii,\,\xi\in{\frak D}$. This form is associated with the Dirichlet
  operator
\begin{equation}\label{z4.4}H_\tsigma^\XR:=H_\tsigma^X+H_\tsigma^\Rp\end{equation}
 on ${\frak D}$ which
  satisfies
\begin{equation}\label{z4.5} {\cal
  E}^\XR_\tsigma(\fii,\xi)=(H^\XR_\tsigma\fii,\xi)
_{L^2(\tsigma)},\qquad
  \fii,\,\xi\in{\frak D}.\end{equation}
Here, $H^X_\tsigma$ and $H^\Rp_\tsigma$ are the Dirichlet operators of
$\nabla^X$ and $\nabla^\Rp$, respectively. They are given by the formulas
\begin{align}
H^X_\tsigma\fii(x,s):&=-\Delta^X\fii(x,s)-\la\nabla^X\log
q(x,s),\nabla^X\fii(x,s)\ra_{T_x(X)},\notag\\
H^\Rp_\tsigma\fii(x,s):&=-\Delta^\Rp\fii(x,s)-\nabla^\Rp\log
q(x,s)\nabla^\Rp\fii(x,s),\label{4.6}\end{align}
where
\begin{align*}\Delta^\Rp&=\div^{\Rp}\nabla^{\Rp}=-(\nabla^{\Rp})^*_{\lambda}\nabla^{\Rp}
\\&=s^2\frac{\di^2}{\di s^2}+2s\frac\di{\di s}=
s^2\frac{\di^2}{\di s^2}+2\nabla^\Rp.\end{align*}

The closure of the form ${\cal E}^\XR_\tsigma$ on $L^2(\pi_\tsigma)$ is denoted
by $({\cal E}^\XR_\tsigma, D({\cal E}^\XR_\tsigma))$. This form generates a
positive selfadjoint operator in $L^2(\tsigma)$ (the so-called Friedrichs
extension of $H_\tsigma^\XR$, see e.g.\ \cite{BKR}). For this extension we
preserve the notation $\HXR$ and denote the domain by $D(\HXR)$.

Using the underlying Dirichlet operator, we obtain the representation
\begin{align}
(\HOM F)(\om):&=-\sum_{j,k=1}^N\frac{\di^2 F}{\di r_j\di r_k
  }(\la\fii_1,\om\ra,\dots,\la\fii_n,\om\ra)
\la\la\nabla^\XR\fii_j,\nabla^\XR\fii_k\ra_{T(\XR)},\om\ra \notag\\
&\quad+\sum_{j=1}^N\frac{\di F}{\di r_j}(\la\fii_1,\om\ra,\dots,\la\fii_n,\om\ra)\la\HXR\fii_j,\om\ra.
\label{4.7}\end{align}

The following theorem implies that  $\HOM$  is well
defined as a linear operator on $\FC$, i.e., independently of the representation
of $F$ as in \eqref{5}.

\begin{th} \label{th4.1}The operator $\HOM$ is associated with the intrinsic Dirichlet form
  $\EOM$\rom, i\rom.e\rom{.,} for all $F,G\in\FC$
\begin{equation}\label{fr}\EOM(F,G)=(\HOM F,G)_{L^2(\pi_\tsigma)},\end{equation}
or
\[ \HOM=-\div^\Om_{\pi_\tsigma}\nabla^\Om\quad {\mathrm on}\ \FC.\]
We call $\HOM$ the intrinsic Dirichlet operator of the measure
$\pi_\tsigma$\rom. \end{th}

\noindent {\it Proof}. For the shortness of notations we will prove the formula
\eqref{fr} in the case where $F,G\in\FC$ are of the form
\[ F=g_F(\la\fii,\om\ra),\quad G=g_G(\la\xi,\om\ra).\]
However, it is a trivial step to generalize the proof for general $F,G$.

Let $\Lambda\in\Oc(X)$ be chosen so that the supports of the functions $\fii$ and $\xi$ are
in $\Lext$. Then, by \eqref{4.1} and \eqref{sex}
\begin{gather*}
\EOM(F,G)=\int_{\OM}g_F'(\la\fii,\om\ra)g'_G(\la\xi,\om\ra)
\la\la\nabla^\XR\fii,\nabla^\XR\xi\ra_{T(\XR)},\om\ra \,\pi_\tsigma(d\om)\\
=-e^{\tsigma(\Lext)}\sum_{n=1}^\infty \frac 1{n!}\int_{\Lambda^n_{\mathrm
    \ext}}g'_F(\fii(x_1,s_1)+\cdots+\fii(x_n,
s_n))g_G'(\xi(x_1,s_1)+\cdots+\xi(x_n,s_n))\\
\times\bigg[\sum_{i=1}^n\la\nabla^{\XR}\fii(x_i,s_i)\nabla^\XR\xi(x_i,s_i)\ra_{T_{(x_i,s_i)}(\XR)}\bigg]
\tsigma(dx_1,ds_1)\dotsm\tsigma(dx_1,ds_1)\\
=e^{-\tsigma(\Lext)}\sum_{n=1}^\infty \frac 1{n!}\int_{\Lambda^n_{\mathrm
    \ext}}\sum_{i=1}^n\la\nabla_i^\XR g_F(\fii(x_1,s_1)+\cdots+\fii(x_n,s_n)),\\
\nabla_i^\XR
g_G(\xi(x_1,s_1)+\dots+\xi(x_n,s_n))\ra_{
  T_{(x_i,s_i)}(\XR)}\,\tsigma(dx_1,ds_1)\dotsm\tsigma(dx_n,ds_n), \end{gather*}
where $\nabla_i^\XR$ denotes the $\nabla^\XR$ gradient in the $(x_i,s_i)$
variables. We note that the vector field
\[\nabla^X_i
g_F(\fii(x_1,s_1)+\cdots+\fii(x_n,s_n))=g'_F(\fii(x_1,s_1)+\cdots+\fii(x_n,s_n))\nabla^X\fii(x_i,s_i)\]
has support in the $x_i$ variable in $\Lambda$ and the vector
field (function)
\[ \nabla_i^\Rp
g_F(\fii(x_1,s_1)+\cdots+\fii(x_n,s_n))=g'_F(\fii(x_1,s_1)+\cdots+\fii(x_n,s_n))
s_i\frac{\di }{\di s_i}\fii(x_i,s_i)\]
is bounded, while $s_iq(x_i,s_i)$ is equal to zero as $s_i=0$ and $s_i=\infty$
for $m$-a.a.\ $x\in X$. Therefore,
\allowdisplaybreaks
\begin{gather*}
\EOM(F,G)=e^{-\tsigma(\Lext)}\sum_{n=1}^\infty\frac
1{n!}\int_{\Lambda^n_{\mathrm \ext}}\bigg[\sum_{i=1}^n
H_\tsigma^{(\XR)_i}g_F(\fii(x_1,s_1)+\cdots+\fii(x_n,s_n))\bigg] \times\\
\times
g_G(\xi(x_1,s_1)+\cdots+\xi(x_n,s_n))\,\tsigma(x_1,s_1)\dotsm\tsigma(dx_n,ds_n)\\
=-e^{\tsigma(\Lext)}\sum_{n=1}^\infty \frac 1{n!}\int_{\Lambda^n_{\mathrm
    \ext}}\bigg[ \sum_{i=1}^n g_F''(\fii(x_1,s_1)+\cdots+\fii(x_n,s_n))\times \\
\times \la
\nabla^\XR\fii(x_i,s_i),\nabla^\XR\fii(x_i,s_i)\ra_{T_{(x_i,s_i)}(\XR)}\\
\mbox{}+g'_F(\fii(x_1,s_1)+\cdots+\fii(x_n,s_n))\big\{\la \nabla^\XR\log
q(x_i,s_i),
\nabla^\XR\fii(x_i,s_i))
\ra_{T_{(x_i,s_i)}(\XR)}\\
\mbox{}+\Delta^X\fii(x_i,s_i)+\Delta^\Rp\fii(x_i,s_i)
\big\}\bigg]\times\\
\times
g_G(\xi(x_1,s_1)+\cdots+\xi(x_n,s_n))\,\tsigma(dx_1,ds_1)\dotsm\tsigma(dx_n,ds_n)\\
=\int_{\OM}\HOM F(\om)G(\om)\,\pi_\tsigma(d\om).\quad\blacksquare
\end{gather*}

\begin{rem}\label{rem4.2}\rom{
The operator $\HOM$ can be naturally extended to cylinder functions of the form
\[ F(\om):=e^{\la\fii,\om\ra},\qquad \fii\in{\frak D},\, \om\in\OM,\]
since such $F$ belong to $L^2(\pi_\tsigma)$. We then have
\begin{equation}\label{4.9}\HOM
  e^{\la\fii,\om\ra}=\la\HXR\fii-|\nabla^\XR\fii|^2_{T(\XR)},\om\ra
 \, e^{\la\fii,\om\ra}.\end{equation}
}\end{rem}

As an immediate  consequence of Theorem~\ref{th4.1} we obtain

\begin{cor}\label{cor4.1}
$(\EOM,\FC)$ is closable on $L^2(\pi_\tsigma)$\rom. Its closure
$(\EOM,D(\EOM))$ is associated with a positive definite selfadjoint
operator\rom, the Friedrichs extension of $\HOM$\rom, which we also denote by
$\HOM$ \rom(and its domain by $D(\HOM)$\rom{).}\end{cor}

Clearly, $\nabla^\Om$
 also extends to $D(\EOM)$. We denote this extension by $\nabla^\Om$.

\begin{cor}\label{cor4.2} Let
\begin{equation}\label{4.10}
\begin{gathered}
F(\om):=g_F(\la\fii_1,\om\ra,\dots,\la\fii_N,\om\ra),\qquad \om\in\OM,\\
\fii_1,\dots,\fii_N\in D({\cal E}^\XR_\tsigma),\ g_F\in C^\infty_{\mathrm
  b}(\R^N).\end{gathered}\end{equation}
Then $F\in D(\EOM)$ and
\[ (\nabla^\Om F)(\omega;x)=\sum_{j=1}^N
\frac{\di g_F}{\di
  r_j}(\la\fii_1,\om\ra,\dots,\la\fii_N,\om\ra)\nabla^\XR\fii_j(x,s_x).\]
\end{cor}

\noindent{\it Proof}. By approximation this is an immediate  consequence of
\eqref{gradient} and the fact that, for all $1\le i\le N$,
\begin{equation}\label{4.11}
\int\la |\nabla^\XR\fii_i|^2_{T(\XR)},\om\ra\,\pi_\tsigma(d\om)={\cal
  E}^\XR_\tsigma(\fii_i,\fii_i).\end{equation}
\begin{rem}\label{rem4.3}\rom{
Let $\mu_{\nu,\tsigma}\in{\cal M}_2(\OM)$ be given as in \eqref{salko}. Then, by
Theorem~3.2, (ii)$\Rightarrow$(i), all results above are
valid with $\mu_{\nu,\tsigma}$ replacing $\pi_\tsigma$. By \eqref{4.7} we have
\[ \HOM=H^\Om_{\mu_{\nu,\tsigma}}\quad\text{on }\FC.\]
We note that the r.h.s.\ of \eqref{4.7} only depends on $\tsigma$ and the
Riemannian structure of $\XR$. The respective Friedrichs extensions on
$L^2(\mu_{\nu,\tsigma})$ again denoted by $H^\Om_{\mu_{\nu,\tsigma}}$, however
do not coincide.
}\end{rem}

\subsection{The heat semigroup and ergodicity}

The results of this subsection are obtained absolutely analogously to the
corresponding results of the paper \cite{AKR}, so we omit the proofs.

\newcommand{\munu}{\mu_{\nu,\tsigma}}
\newcommand{\HH}{H^\XR_{\tsigma}}
\newcommand{\TTT}{T^\Om_{\munu}(t)}
\newcommand{\HHH}{H^\Om_{\munu}}
\newcommand{\EE}{{\cal E}_\tsigma^\XR}
\newcommand{\EEE}{{\cal E}^\Om_{\munu}}

For $\munu\in {\cal M}_2(\OM)$ let $\TTT:=\exp(-t\HHH)$, $t>0$. Define
\[ E({\frak
  D}_1,\OM)=\operatorname{l.h.}\big\{\,\exp(\la\log(1+\fii),\cdot\ra)\mid\fii\in{\frak
  D}_1\,\big\},\]
where l.h.\ means the linear hull and
\begin{align*}
{\frak D}_1:=\big\{\,\fii\in D(\HH)&\cap L^1(\tsigma)\mid \HH\fii\in
L^1(\tsigma)\\
&\text{and }-\delta\le\fii\le0\text{ for some }\delta\in(0,1)\,\big\}
.\end{align*}

\begin{prop}\label{prop4.1}
Let $\munu$ be as in \eqref{salko}\rom. Assume that $\HH$ is conservative\rom,
i\rom.e\rom{.,}
\[\int_{\XR}(\HH\fii)(x,s)\,\tsigma(dx,ds)=0\]
for all $\fii\in D(\HH)\cap L^1(\tsigma)$ such that $\HH\fii\in
L^1(\tsigma)$\rom, and suppose that $(\HH,{\frak D})$ is essentially selfadjoint
on $L^2(\tsigma)$\rom. Then
\begin{equation}\label{peter}
\TTT\exp(\la\log(1+\fii),\cdot\ra)=\exp(\la\log(1+e^{-t\HH}\fii),\cdot\ra),\qquad
\fii\in{\frak D}_1,\end{equation}
$E({\frak D}_1,\OM)\subset D(\HHH)$\rom, and
\begin{align*}
&\HHH\exp(\la\log(1+\fii),\cdot\ra)\\
&\qquad =\la (1+\fii)^{-1}\HH\fii,\cdot\ra\exp(\la\log(1+\fii),\cdot\ra),\qquad
\fii\in{\frak D}_1.\end{align*}
\end{prop}

\begin{rem}\rom{
(i) The condition of selfadjointness of $\HH$ on ${\frak D}$  is fulfilled if
$X$ is complete and $|\beta^\tsigma|_{T(\XR)}\in L^p_{\mathrm loc}(X;m)\otimes
L^p(\Rp;\lambda)$ for some $p\ge\operatorname{dim}(X)+1$.

(ii) Since $(\exp(-t \HH))_{t>0}$ is sub-Markovian (i.e.,
$0\le\exp(-t\HH)\fii\le1$ for all $t>0$ and $\fii\in L^2(\tsigma)$,
$0\le\fii\le1$), because $(\EE,D(\EE))$ is a Dirichlet form, by a simple
approximation argument Proposition~\ref{prop4.1}  implies that the equality
\eqref{peter} holds for $t>0$ and
all $\fii\in L^1(\tsigma)$, $-1<\fii\le0$.
}\end{rem}

\begin{th}\label{th4.2}
Let the conditions of Proposition~\rom{\ref{prop4.1}} hold\rom. Then $E({\frak
  D}_1,\OM)$ is an operator core for the Friedrichs extension $\HHH$ on
  $L^2(\munu)$\rom. \rom(In other words\rom:\linebreak $(\HHH,E({\frak D}_1,\OM))$ is
  essentially selfadjoint on $L^2(\munu)$\rom{.)}
\end{th}

\begin{th}\label{th4.3} Suppose that the conditions of
  Theorem~\rom{\ref{hahaha}}  and Proposition~\rom{\ref{prop4.1}} hold\rom. Then the
  following assertions are equivalent\rom:

\rom{(i)} $\munu=\pi_{z\tsigma}$ for some $z>0$\rom.

\rom{(ii)} $(\EEE,D(\EEE))$ is irreducible \rom(i\rom.e\rom{.,} for $F\in
D(\EEE)$\rom, $\EEE(F,F)=0$ implies that $F=\operatorname{const}$\rom{).}

\rom{(iii)} $(\TTT)_{t>0}$ is irreducible \rom(i\rom.e\rom{.,} if $G\in
L^2(\munu)$ such that $\TTT(GF)=G\TTT F$ for all $F\in L^\infty(\munu)$\rom,
$t>0$\rom, then $G=\operatorname{const}$\rom{).}

\rom{(iv)} If $F\in L^2(\munu)$ such that $\TTT F=F$ for all $T>0$\rom, then $F=\operatorname{const}.$

\rom{(v)} $\TTT\not\equiv\pmb{\pmb{1}}$ and ergodic \rom(i\rom.e\rom{.,}
\[\int\bigg(\TTT F-\int F\,d\munu\bigg)^2d\munu\to0\quad \mathrm{as}\ t\to0\]
for all $F\in L^2(\munu)$\rom{).}

\rom{(vi)} If $F\in D(\HHH)$ with $\HHH=0$\rom, then $F=\operatorname{const}.$
\end{th}

\begin{rem}\rom{
Let us consider the diffusion process $M$ on $\XR$ associated to the Dirichlet
form $(\EE,D(\EE))$. This process can be interpreted as distorted Brownian
motion on the manifold $\XR$. More precisely, the diffusion of points $x\in X$
is associated to the Dirichlet form of the measure $\sigma$,
 so that it is distorted Brownian motion on
$X$. The diffusion of marks $s_x$, $x\in X$, can be obtained by the exponential
change of time $\R\ni t\mapsto {\cal L}^{-1}t=e^t=s\in\Rp$ in the distorted
Brownian motion on $\R$ associated to the Dirichlet form of the measure ${\cal
  L}^*p(x,\cdot)$. This follows from the representation \eqref{4.5}--\eqref{4.6}
of the Dirichlet form $\EE$ and Remark~\ref{rem3.1}.

The existence of a diffusion process ${\bf M}$ corresponding to the Dirichlet
form $(\EEE,D(\EEE))$ and its identification with the independent infinite
particle process (on $\XR$)
may be proved by the same arguments as in \cite{AKR}. By
analogy with the case of the process $M$ on $\XR$, one can call $\bf M$
distorted Brownian motion on $\OM$.
}\end{rem}

\section{Intrinsic Dirichlet operator and second quantization}

\renewcommand{\TTT}{T^\Om_{\pi_\tsigma}(t)}
\renewcommand{\HHH}{H^\Om_{\pi_\tsigma}}
\renewcommand{\EEE}{{\cal E}^\Om_{\pi_\tsigma}}
\newcommand{\nMP}{\nabla^{\mathrm MP}}
\newcommand{\ot}{\otimes}
\newcommand{\hot}{\widehat\otimes}

In this section, we want to describe the Fock space realization of the marked
Poisson spaces and show that $\HHH$ is the second quantization of the operator
$\HH$.

\subsection{Marked Poisson gradient and chaos decomposition}

Let us define another ``gradient'' on functions $F\colon \OM\to\R$, which has
specific useful properties on the marked Poisson space.

\begin{define}\rom{
For any $F\in\FCp$ we define the marked Poisson gradient
$\nMP$ as
\[ (\nMP F)(\om,(x,s)):=F(\om+\eps_{(x,s)})-F(\om),\qquad \om\in\OM,\ (x,s)\in \XR.\]
}\end{define}

Let us mention that the operation
\[\OM\ni\om\mapsto\om+\eps_{(x,s)}\in\OM\]
is a $\pi_\tsigma$-a.e.\ well-defined map because of the property
\[ \pi_\tsigma\big(
\{\om=(\gamma,s)\in\OM\mid x\in\gamma\}
\big)=0\]
for an arbitrary $x\in X$ (which easily follows from the construction of
$\pi_\tsigma$). We consider $\nMP$ as a mapping
\[\nMP\colon \FCp\ni F\mapsto \nMP F\in L^2(\tsigma)\otimes L^2(\pi_\tsigma)\]
that corresponds to using the Hilbert space $L^2(\tsigma)$ as a tangent space at
any point $\om\in\OM$. Thus, for any $\fii\in{\frak D}$, we can introduce the
directional derivative
\begin{align*}
(\nabla^{\mathrm MP}_\fii F)(\om)&=\la\nMP F(\om),\fii\ra_{L^2(\tsigma)}\\
&=\int_{\XR}(F(\om+\eps_{(x,s)})-F(\om))\fii(x,s)\,\tsigma(dx,ds)
.\end{align*}

The most important feature of the marked Poisson gradient is that it produces
(via a corresponding ``integration by parts formula'') the orthogonal system of
Charlier polynomials on $(\OM,{\cal B}(\OM),\pi_\tsigma)$. Below, we describe this
construction in detail using the isomorphism between $L^2(\pi_\tsigma)$ and the
Fock space (see \cite{ItKu,KSS2,LRS})

Let ${\cal F}(L^2(\tsigma))$ denote the symmetric Fock space over
$L^2(\tsigma)$:
\[ {\cal F}(L^2(\tsigma)):=\bigoplus_{n=0}^\infty {\cal F}_n(L^2(\tsigma))n!,\]
where
\[ {\cal F}_n(L^2(\tsigma)):=(L^2(\tsigma))^{\hot n}=\hat
L^2((\XR)^n,\tsigma^{\ot n})\]
and ${\cal F}_0(L^2(\tsigma)):=\R$.
Thus, for each $F=(f^{(n)})_{n=0}^\infty\in {\cal F}(L^2(\tsigma))$
\[ \| F\|^2_{{\cal F}(L^2(\tsigma))}=\sum_{n=0}^\infty |f^{(n)}|_{\hat
  L^2(\tsigma^{\ot n})}n!.\]

By ${\cal F}_{\mathrm fin}({\frak D})$ we denote the dense subset of ${\cal
  F}(L^2(\tsigma))$ consisting of finite sequences $(f^{(n)})_{n=0}^N$,
  $n\in\Z_+$, such that each $f^{(n)}$ belongs to  ${\cal F}_n({\frak D}):={\mathrm
  a}.{\frak D}^{\hot n}$, the $n$-th symmetric algebraic tensor power of ${\frak
  D}$:
\[ {\mathrm a}.{\frak D}^{\hot n}:=\operatorname{l.h.}\{\fii_1\hot\dotsm\hot
  \fii_n\mid\fii_i\in{\frak D}\}.\]
In virtue of the polarization identity, the latter  set is spanned just by the
vectors of the form $\fii^{\ot n}$ with $\fii\in{\frak D}$.

Now, we define a linear mapping
\begin{equation}{\cal F}_{\mathrm fin}({\frak D})\ni F=(f^{(n)})_{n=0}^N\mapsto
IF=(IF)(\om)=\sum_{n=0}^N Q_n (f^{(n)};\om)\in\FP
\label{mapping}\end{equation}
by using the following recursion relation:
\begin{align}
Q_{n+1}(\fii^{\ot(n+1)};\om)&=Q_n(\fii^{\ot
  n};\om)(\la\om,\fii\ra-\la\fii\ra_{\tsigma})\notag\\ &\quad-
nQ_n(\fii^{\ot(n-1)}\hot(\fii^2),\om)-nQ_{n-1}(\fii^{\ot(n-1)};\om)\la\fii^2\ra_{\tsigma},\notag\\
&\quad Q_0(1,\om)=1,\quad \fii\in{\frak D}.\label{ser}
\end{align}
Here, we denoted by $\la\fii\ra_{\tsigma}:=\int\fii\,d\tsigma$. Notice that, since ${\frak
  D}$ is an algebra under pointwise multiplication of functions, the latter
  definition is correct.

  It is not hard to see that the mapping \eqref{mapping} is
  one-to-one. Moreover, the following proposition holds:

\begin{prop}\label{prop5.1} The mapping \eqref{mapping} can be extended by
  continuity to a unitary isomorphism between the spaces ${\cal
  F}(L^2(\tsigma))$ and $L^2(\pi_\tsigma).$\end{prop}

For each $\fii\in{\frak D}$, let us define the creation and annihilation
operators in ${\cal F}(L^2(\tsigma))$ by
\[ a^+(\fii)\psi^{\ot n}=\fii\hot\psi^{\ot n},\qquad a^-(\fii)\psi^{\ot
  n}=n(\fii,\psi)_{L^2(\tsigma )}\psi^{\ot (n-1)},\qquad \psi\in{\frak D}.\]

We will denote by the same letters the images of these operators under the
unitary $I$.

\begin{prop}\label{prop5.2} We have\rom, for each $\fii\in{\frak D},$
\[ a^-(\fii)=\nabla^{\mathrm MP}_\fii,\qquad a^+(\fii)=\nabla_\fii^{{\mathrm
    MP}*}.\]
In particular\rom,
\[ Q_n(\fii_1\hot\dotsm\hot\fii_n;\om)=(\nabla^{{\mathrm MP}*}_{\fii_1}\dotsm
\nabla^{{\mathrm MP}*}_{\fii_n}1)(\om),\qquad \om\in\OM.\]
\end{prop}

\newcommand{\Exp}{\operatorname{Exp}}

Finally, for each $\fii\in{\frak D}$ we introduce the Poisson exponential
\[ e(\fii;\cdot):=\sum_{n=0}^\infty \frac 1{n!}\, Q_n(\fii^{\ot
  n};\cdot)=I(\Exp\fii),\]
where
$$\Exp\fii=\bigg(\frac1{n!}\,\fii^{\ot n}\bigg)_{n=0}^\infty.$$
Then, one can  show that, for $\fii>-1$,
\begin{equation}\label{zebra}
e(\fii;\om)=\exp\big[\la\log(1+\fii),\om\ra-\la\fii\ra_{\tsigma}\big],\qquad
\om\in\OM. \end{equation}

\subsection{Second quantization on the marked Poisson space}

Let $B$ be a contraction on $L^2(\tsigma)$, i.e., $B\in {\cal
  L}(L^2(\tsigma),L^2(\tsigma))$, $\|B\|\le1$. Then, we can define the operator
  $\Exp B$ as the contraction on ${\cal F}(L^2(\tsigma))$ given by
\begin{align*}
&\Exp B \restriction {\cal F}_n(L^2(\tsigma)):=B\ot\dotsm\ot B\quad \text{($n $
  times)},\ n\in\N,\\
&\Exp B\restriction {\cal F}_0(L^2(\tsigma)):=1.\end{align*}

For any selfadjoint positive operator $A$ in $L^2(\tsigma)$, we have a contraction
semigroup $e^{-tA}$, $t\ge0$, and it is possible to introduce a positive
selafadjoint operator $d \Exp A$ as the generator of the semigroup
$\Exp(e^{-tA})$, $t\ge0$:
\begin{equation}\label{kon}
\Exp(e^{-tA})=\exp(-td\Exp A).\end{equation}
The operator $d\Exp A$ is called the second quantization of $A$. We denote by
$H^{\mathrm MP}_A$ the image of the operator $d\Exp A$ in the marked Poisson
space $L^2(\pi_\tsigma)$.

\begin{th}\label{th5.1}
Let ${\frak D}\subset\operatorname{Dom}A$\rom. Then\rom, the symmetric bilinear
form corresponding to the operator $H_A^{\mathrm MP}$ has the following
representation\rom:
\begin{equation}\label{5.18}
(H_A^{\mathrm MP}F,G)_{L^2(\pi_\tsigma)}=\int_{\OM}(\nMP F,A\nMP
G)_{L^2(\tsigma)}\,\pi_\tsigma(d\om) \end{equation}
for all $F,G\in\FP.$\end{th}

\begin{rem}\label{rem5.1}
\rom{The bilinear form \eqref{5.18} uses the marked Poisson gradient $\nMP$ and
  a coefficient operator $A>0$. We will call
\[ {\cal E}^{\mathrm MP}_{\pi_\tsigma,A}(F,G)=\int_{\OM}(\nMP F,A\nMP
  G)_{L^2(\tsigma)}\,\pi _{\tsigma}(d\om)\]
the marked Poisson pre-Dirichlet form with coefficient $A$.}\end{rem}

\noindent {\it Proof of Theorem\/} 5.1. The proof is analogous to that of
Theorem~5.1 in \cite{AKR}. Using again the fact that ${\frak D}$ is an algebra
under pointwise multiplication, one easily concludes that, for any $F\in\FP$ and
any $\om\in\OM$, the gradient $\nMP F(\om,(x,s))$ is
a function in $\frak D$ and hence
\[ (\nMP F,A\nMP G)_{L^2(\tsigma)}\in\FP,\]
so that the form \eqref{5.18} is well-defined. Then, one verifies the formula
\eqref{5.18} by using Propositions~5.1, 5.2 and the explicit formula
for $d\Exp A$ on ${\cal F}_n({\frak D})$:
\[ d\Exp A\,\fii^{\ot n}=n(A\fii)\hot \fii^{\ot(n-1)},\qquad \fii\in{\frak
  D}.\quad
\blacksquare\]

\subsection{The intrinsic Dirichlet operator as a second quantization}

The following two theorems are again analogous to corresponding results
(Theorems~5.2 and~5.3) in \cite{AKR}, so we omit their proofs.

Let us consider the special case of the second quantization operator $d\Exp A$
where  the operator $A$ coincides with the Dirichlet operator $H^{\XR}_\tsigma$.

\begin{th}\label{th5.2} We have the equality
\[ H^{\mathrm MP}_{\HH}=\HHH\]
on the dense domain $\FCp$\rom. In particular\rom, for all $F,G\in\FCp$
\begin{align*}
& \int_{\OM}\la \nabla^\Om F(\om),\nabla^\Om
G(\om)\ra_{T_\om(\OM)}\,\pi_\tsigma(d\om)\\
&\qquad = \int_{\OM}(\nMP F(\om),\HH\nMP
G(\om))_{L^2(\tsigma)}\,\pi_\tsigma(d\om),\end{align*}
or
\[ \nabla^{\Om*}\nabla^\Om=\nabla^{{\mathrm MP}*}\HH\nMP\]
as an equality on $\FCp.$\end{th}

\begin{th}\label{th5.3}Suppose that the operator $\HH$ is essentially
  selfadjoint on the domain ${\frak
  D}\subset\operatorname{Dom}(H_\tsigma^\XR)$\rom. Then\rom, the intrinsic
  Dirichlet operator $\HHH$ is essentially selafadjoint on the domain
  $\FC.$\end{th}

\begin{rem}\label{rem5.2}\rom{
Notice that in Theorem~\ref{th5.3} we do not suppose the operator $\HH$ to be
conservative. So, this theorem is a generalization of Theorem~\ref{th4.2} in the
special case where $\munu=\pi_\tsigma$.}\end{rem}

\begin{cor}\label{cor5.1}
Suppose that the condition of Theorem~\rom{\ref{th5.3}}
is satisfied and let $T^\Om_{\pi_\tsigma}(t)=\exp(-t\HHH)$\rom,
$t>0$\rom. Then\rom, for each $\fii\in{\frak D}$\rom, $\fii>-1$\rom, we have
\begin{equation}
\TTT\exp(\la\log(1+\fii),\cdot\ra)
=
\exp\big[\la\log(1+e^{-t\HH}\fii),\cdot\ra-\la(e^{-t\HH}-\pmb{1})\fii\ra_\tsigma\big].
\label{mein}\end{equation}\end{cor}

\noindent {\it Proof}. The formula \eqref{mein} follows from
Proposition~\ref{prop5.1}, \eqref{zebra}, \eqref{kon} and Theorems~\ref{th5.2}
and~\ref{th5.3}.\quad$\blacksquare$

\begin{rem}\rom{
If $\HH$ is conservative, then
\[ \int(e^{-t\HH}-\pmb{1})\fii\,\,d\tsigma=0\qquad\text{for all }t\ge0,\]
and so in this case \eqref{mein} coincides with \eqref{peter}
for $\fii\in{\frak D}$, $\fii>-1$.}\end{rem}

\section{Appendix}

We will prove now the part (i)$\Rightarrow$(ii) of Theorem~\ref{hahaha}. To this end, we
present first some definitions and reformulation of known statements.

For any $\Upsilon\in{\cal B}(X)$ and $\omega\in\OM$, denote
$\om_\Upsilon:=\om\cap\Upsilon$.

\begin{define}\rom{
For $\Lambda\in\Oc(X)$, we define for $\om\in\OM$ and $\Delta\in{\cal B}(\OM)$
\begin{equation}\label{a1}
\widehat\Pi^\tsigma_{\Lext}(\om,\Delta):=\int_{\OM}\pmb{\pmb{1}}_\Delta(\om_{(X\setminus\Lambda)_{\mathrm
    \ext}}+\om'_{\Lext})\,\pi_{\tsigma}(d\om'\mid
    N_{\Lext}(\cdot)=\om(\Lext)).\end{equation}
A probability measure $\mu$ on $(\OM,{\cal B}(\OM))$ is called a canonical Gibbs
    measure for the free case if
\[\mu(\Delta)=\int_{\OM}\widehat\Pi^\tsigma_{\Lext}(\om,\Delta)\,\mu(d\om).\]
Let ${\cal G}_{\mathrm c}(\tsigma)$ denote the set of all such probability measures.
}\end{define}

For $\Lambda\in\Oc(X)$ denote by
$\Bhat_{(X\setminus\Lambda)_{\mathrm \ext}}(\OM)$ the $\sigma$-algebra generated
by the $\sigma$-algebra ${\cal B}_{(X\setminus\Lambda)_{\mathrm \ext}}(\OM)$ and
the mappings $N_{\Lext}$. Then, one easily deduces the following proposition.

\begin{prop}\label{A1} A probability measure $\mu$ on $(\OM,{\cal B}(\OM))$
  belongs to ${\cal G}_{\mathrm c}(\tsigma)$ if and only if for each bounded
  measurable function $G$ on $\OM$
\[
{\Bbb E}_\mu[G\mid\Bhat_{(X\setminus\Lambda)_{\mathrm \ext}}(\OM)]=\widehat
\Pi^\tsigma_{\Lext}G\quad \text{\rom{$\mu$-a.e.}}
\]
Here\rom, for a sub-$\sigma$-algebra $\Sigma\subset{\cal B}(\OM)$\rom, ${\Bbb
  E}_\mu[\cdot \mid\Sigma]$ denotes the conditional expectation
  w\rom.r\rom.t\rom.\ $\mu$ and $\Sigma$ and
\[
  \widehat\Pi^\tsigma_{\Lext}G=(\widehat\Pi^\tsigma_{\Lext}G)(\om):=\int_{\OM}G(\om')\, \widehat\Pi
^\tsigma_{\Lext}(\om,d\om').\]
\end{prop}

The following theorem, which is in fact due to \cite{NZ}, see also
\cite{Georgii} and \cite{KMM}, can be obtained from the original one by a simple
modification of the proof.

\begin{th}\label{A2} Let $\mu$ be a probability measure on $(\OM,{\cal
    B}(\OM))$\rom. Then\rom, $\mu\in{\cal G}_{\mathrm c}(\tsigma)$ if and only if
    there exists a probability measure $\nu$ on $(\R_+,{\cal B}(\R_+))$ such
    that
\[\mu=\int_{\R_+}\pi_{z\tsigma}\,\nu(dz).\]
\end{th}

Hence, by virtue of Proposition~\ref{A1} and Theorem~\ref{A2}, it suffices to
show that if $\mu\in({\mathrm IbP})^\tsigma$ then for $\mu$-a.e.\ $\om\in\OM$
\begin{multline}\label{a4}
{\Bbb E}_\mu [F\mid \Bhat_{(X\setminus\Lambda)_{\mathrm
    \ext}}](\om)=\tsigma(\Lext)^{-\om(\Lext)}\times\\
\times
\int_{\Lext}\dotsi\int_{\Lext}
F(\eps_{(x_1,s_1)}+\dots+
\eps_{(x_{\om(\Lext)},s_{\om(\Lext)})})\,\tsigma(dx_1,ds_1)\dots
\tsigma(dx_{\om(\Lext)},ds_{\om(\Lext)})\end{multline}
for all bounded ${\cal B}_{\Lext}(\OM)$ measurable functions $F\colon\OM\to\Rp$
and all $\Lambda\in\Oc(X)$. Here, the r.h.s.\ of \eqref{a4} is understood to be
equal to $F(\varnothing)$ if $\om(\Lext)=0$. So, fix $\Lambda\in\Oc(X)$.

{\it Claim}. Let $g=(v,a)\in V_0(\Lambda)\times
C_0(\Lambda)$, let $K_1$,\dots,$K_M$ be arbitrary subsets of
$(X\setminus\Lambda)_{\mathrm \ext}$ from $\Bc(\XR)$, and let
$\fii_1,\dots,\fii_N\in C^\infty_{\text{0,b}}(\Lext)$. Let
$G:=g_G(N_{K_1},\dots,N_{K_M})$,
$F:=g_F(\la\fii_1,\cdot\ra,\dots,\la\fii_N,\cdot\ra)$ with $g_G\in C_{\mathrm
  b}^\infty(\R^N)$, $g_F\in C^\infty_{\mathrm b}(\R^M)$ and $n\in\N$. Then
\[\int G{\pmb 1}_{\{N_{\Lext}=n\}}\nabla^\Omega_{(v,a)}F\,f\mu=-\int G{\pmb
  1}_{\{N_{\Lext}=n\}}
B^{\pi_\tsigma}_{(v,a)}F\,d\mu.\]

To prove the claim, for $l\in\N$, $1\le m\le M$, we choose functions
$\chi_l,\chi_{ml}\in C^\infty_{\mathrm 0,b}(\XR)$ taking values in $[0,1]$ such
that

a) $\chi_l=1$ on a neighborhood of $K_g$ (see Section~2), $\chi_l\le {\pmb
  1}_{\Lext}$, and $\chi_l\to{\pmb 1}_{\Lext}$ as $l\to\infty$.

b) $\chi_{ml}=1$ on $K_m$, $\chi_{ml}=0$ in a neighborhood of $K_g$, and
$\chi_{ml}\to{\pmb 1}_{K_m}$ as $l\to\infty$.

Furthermore, let $g\in C_0^\infty(\R)$ be such that ${\pmb 1}_{\{n\}}\le
g\le{\pmb 1}_{]n-\frac12,n+\frac12[}$. Then, for every $\om\in\OM$
\[ G_l(\om):=g_G(\la \chi_{1l},\om\ra,\dots,\la\chi_{Ml},\om\ra)=G(\om)\]
and
\[ g_l(\om):=g(\la\chi_l,\om\ra)={\pmb 1}_{\{N_{\Lext}=n\}}(\om)\]
for all $l>l(\om)$. Moreover, for all $l\in\N$, $1\le m\le M$, we have that
$\nabla^{\XR}_{(v,a)} \chi_l\equiv 0\equiv
\nabla^{\XR}_{(v,a)}\chi_{ml}$. Hence,
\[ \nabla^\Omega_{(v,a)}g_l=0=\nabla^\Omega_{(v,a)}G_l.\]
Consequently, since $\mu\in(\text{IbP})^\tsigma$,
\begin{gather*}
\int G{\pmb 1}_{\{N_{\Lext}=n\}}\nabla^\Omega_{(v,a)}F\,d\mu=\lim_{l\to\infty}\int
G_lg_l\nabla^\Omega_{(v,a)}F\, d\mu\\
=-\lim_{l\to\infty}\bigg[
\int(g_l\nabla^\Omega_{(v,a)}G_l+G_l\nabla^\Omega_{(v,a)}g_l)F\,d\mu+
\int G_lg_lFB^{\pi_\tsigma}_{(v,a)}\,d\mu\bigg]\\
=-\int G{\pmb 1}_{\{N_{\Lext}=n\}}FB^{\pi_\tsigma}_{(v,a)}\,d\nu,\end{gather*}
and the claim is proven.

The claim immediately implies that for $F$ and $(v,a)$ as in the claim
\begin{equation}\label{a5}
{\Bbb E}_\mu[\nabla^\Omega_{(v,a)}F\mid \Bhat_{(X\setminus\Lambda)_{\mathrm
    \ext}}(\OM)]= -{\Bbb E}_\mu[B^{\pi_\tsigma}_{(v,a)}F\mid
    \Bhat_{(X\setminus\Lambda)_{\mathrm \ext}}(\OM)]\quad \text{$\mu$-a.e.}\end{equation}

Now, set $A_1:=\Lext$ and pick $A_i\in\Bc((X\setminus\Lambda)_{\mathrm \ext})$,
$i\in\N$, $i\ge2$, closed under finite intersections such that
\[ T:=(N_{A_i})_{i\in\N}\colon\OM\to\Z_+^\infty\]
generates $\Bhat_{(X\setminus\Lambda)_{\mathrm \ext}}(\OM)$. Disintegrating $\mu$
w.r.t.\  $T$, there exists a probability kernel $\tilde\mu\colon
\Z_+^\infty\times{\cal B}(\OM)\to[0,1]$ such that
\begin{equation}\label{a6}
\mu=\tilde\mu(\un,d\om)(\mu\circ T^{-1})(d\un),\end{equation}
i.e., $\tilde\mu(\un,d\om)$ is a regular conditional probability corresponding to
${\Bbb E}_\mu[\cdot\mid T=\un]$, $\un\in\Z_+^\infty$. From now on all statements
depending on $\un$ below are meant as statements which are true $(\mu\circ
T^{-1})(d\un)$-a.e. It follows from \eqref{a5} that
\begin{equation}\label{a7}
\int\nabla^\Omega_{(v,a)}F(\om)\,\tilde\mu(\un,d\om)=-\int\la\beta^{\tsigma}_{(v,a)},\om\ra
F(\om)\,\tilde\mu(\un,d\om)\end{equation}
for all $(v,a)$ and $F$ as in the claim. Furthermore, since clearly
$\tilde\mu(\un,d\om)$ is supported by $\{T=\un\}$ and since for $n\in\N$
$\{\om\subset\Lext\mid \om(\Lext)=n\}$ is isomorphic (as a measurable space) to
${\widetilde\Lambda}_{\text{\ext}}^n/{\frak S}_n$, there exists a probability
measure $\mu_{\un}$ on ${\widetilde \Lambda}_{\mathrm \ext}^n$ invariant under
${\frak S}_n$ such that for all positive ${\cal B}(\OM)$-measurable functions
$G$ on $\OM$
\begin{equation}
\label{a8}
\int G(\om)\,\tilde\mu(\un,d\om)=\int_{{\widetilde \Lambda}_{\mathrm
    \ext}^n}G\big( \sum_{k=1}^n \eps
    _{(x_k,s_k)}+{\tilde\om}_{(X\setminus\Lambda)_{\mathrm
    \ext}}\big)\,\mu_{\un}(dx_1,ds_1,\dots,dx_n,ds_n), \end{equation}
where $\un =(n,n_2,n_3,\dots)$ and $\tilde\om_{(X\setminus\Lambda)_{\mathrm
    \ext}}$ is the unique element in $\Omega^{\Rp}_{(X\setminus\Lambda)_{\mathrm
    \ext}}$ such that $\tilde\om_{(X\setminus\Lambda)_{\mathrm
    \ext}}(A_i)=n_i$ for all $i\ge2$. By the invariance under ${\frak S}_n$, it
    suffices to determine $\mu_{\un}$ on just one of the $n!$ connected
    components of ${\widetilde\Lambda}^n_{\mathrm \ext}$. Let
    ${}_{\circ}{\widetilde\Lambda} ^n _{\mathrm \ext}$ denote this component.
Hence, it suffices to determine $\mu_{\un}$ on
    $M:={}_1\Lext\times\dotsm\times{}_n\Lext\subset
    {}_{\circ}{\widetilde\Lambda}^n_{\mathrm \ext}$, ${}_i\Lambda$, $1\le i\le n$,
    pairwise disjoint open subsets of $\Lambda$. Therefore, by
    Lemma~\ref{lemcharacter} and a simple disintegration argument, it is enough
    to show that, for all $\fii_i\in C^\infty_{\mathrm 0,b}({}_i\Lext)$ and all
    $(v_i,a_i)\in V_0({}_i\Lambda)\times C_0 ({}_i\Lambda)$, $1\le i\le n$,
\begin{gather}
\int\fii_1(x_1,s_1)\dotsm\la\nabla^{\XR}\fii_i,(v_i,a_i)\ra_{\XR}(x_i,s_i)\dotsm(x_n,s_n)\,
\mu_{\un}(dx_1,ds_1,\dots,dx_n,ds_n)\notag\\
=-\int\beta^\tsigma_{(v_i,a_i)}\fii_1(x_1,s_1)\dotsm\fii_n(x_n,s_n)\,\mu_{\un}(dx_1,ds_1,\dots,dx_n,ds_n).\label{a9}\end{gather}
Because then by Lemma~\ref{lemcharacter} we see that $\mu_{\un}$ is up to a
constant equal to $\tsigma^n$ on ${}_1\Lext\times\dotsm\times{}_n\Lext$, hence
since $\mu_{\un}$ is a probability measure on ${\widetilde \Lambda}^n_{\mathrm
  \ext}$, it follows that
\[\mu_{\un}=\tsigma(\Lext)^{-n}\tsigma^n,\]
and  \eqref{a4} follows then from \eqref{a6} since the case $n=0$ is trivial. To
show \eqref{a9} we fix $i\in\{1,\dots,n\}$. Choosing $F$ in \eqref{a7} so that
$F:=\prod_{j=1} ^n\la\fii_j,\cdot\ra$ on $\{\om\subset\Lext\mid\om(\Lext)=n\}$,
we conclude that
\begin{gather*}
\int\la\fii_1,\om\ra\dotsm\la\nabla^{\XR}_{(v_i,a_i)}\fii_i,\om\ra\dotsm\la\fii_n,\om\ra\,\tilde\mu(\un,
d\om)\\
=-\int\la\beta^\tsigma_{(v_i,a_i)},\om\ra\la\fii_1,\om\ra\dotsm\la\fii_n,\om\ra\,\tilde\mu(\un,d\om).
\end{gather*}
Hence, by \eqref{a8}

\begin{align}\notag
&\int_{{\widetilde\Lambda}^n_{\mathrm \ext}}
\prod_{\substack{j=1,\\ j\ne
  i}}
(\fii_j(x_1,s_1)+\cdots+\fii_j(x_n,s_n)
  )\sum_{k=1}^n\nabla^{\XR}_{(v_i,a_i)}\fii_i(x_k,s_k)\,
  \mu_{\un}(dx_1,ds_1,\dots,dx_n,ds_n)\\
&\qquad \notag =-\int_{{\widetilde\Lambda}^n_{\mathrm
  \ext}}\sum_{k=1}^n\beta^\tsigma_{(v_i,a_i)}(x_k,s_k)\times\\
&\qquad\qquad\times \prod_{j=1}^n
  (\fii_j(x_1,s_1)+\cdots+\fii_j(x_n,
  s_n))\,\mu_{\un}(dx_1,ds_1,\dots,dx_n,ds_n).
\label{a10}\end{align}
Since both integrals are invariant under ${\frak S}_n$, \eqref{a10} also holds
  if we inly take the integrals over ${}_{\circ}{\widetilde\Lambda}^n_{\mathrm
  \ext}$. But then \eqref{a10} directly turns into \eqref{a9}, since $\fii_i$,
  $a_i$ have support in $\Lambda_i$, $1\le i\le n$, which in turn are pairwise
  disjoint.\quad $\blacksquare$

\begin{center}\bf ACKNOWLEDGMENTS\end{center}

One of the authors, Yu.K., was partially supported through the INTAS-Project
Nr.~97--0378. E.L.  acknowledges financial support of the
Graduiertenkolleg of the University of Bielefeld.
A partial financial support of DFG through the project
436 UKR 113/39/1 is gratefully acknowledged by G.U.

\end{document}